\documentclass[11pt,fullpage]{article}
\usepackage{amsfonts}

\font\bbigsym=cmmi10 scaled\magstep4 
\font\bbbigsym=cmmi10 scaled\magstep5 
\newcommand{\BIGP}{\lower2pt\hbox{\bbigsym\char'031}} 
\newcommand{\BBIGP}{\lower3pt\hbox{\bbbigsym\char'031}}

\newcommand{\openE}{\hbox{I\kern-.2em E}} 
\newcommand{\openN}{\hbox{I\kern-.2em N}} 
\newcommand\openR{\hbox{I\kern-.2em R}} 
\newcommand\smopenR{{\hbox{\sevenrm I \kern-.5em R}}} 
\newcommand\openH{\hbox{I\kern-.2em H}} 
\newcommand\openL{\hbox{I\kern-.2em L}}

\renewcommand\epsilon{\varepsilon} 
\renewcommand \phi{\varphi}

\renewcommand{\hat}{\widehat} 
\renewcommand\tilde{\widetilde} 
\renewcommand{\bar}{\overline}

\newcommand\Verti{\Vert_{\infty}}
\newcommand\Vertv{\Vert_v}

\newcommand{\cov}{{\rm cov \;}}
\newcommand{\var}{{\rm var \;} }

\newcommand{\sign}{{\rm sign \;}}

\newcommand{\diag}{{\rm diag \;}}

\newcommand{\done}{$ \ \Box$ }

\newcommand{\B}{{\cal B}}

\newcommand{\D}{{\cal D}}

\newcommand{\F}{{\cal F}}
\newcommand{\G}{{\cal G}}
\renewcommand{\H}{{\cal H}}
\newcommand{\I}{{\cal I}}

\newcommand{\N}{{\cal N}}
\renewcommand{\P}{{\cal P}}

\newcommand{\T}{{\cal T}}

\newcommand{\V}{{\cal V}}

\renewcommand{\l}{{\ell}}

\def\for{\quad {\rm for} \quad}

\def\tiff{\quad{\rm iff} \quad}

\begingroup

\newtheorem{thm}{Theorem}[section]
\newtheorem{prop}[thm]{Proposition}
\newtheorem{lem}[thm]{Lemma}
\newtheorem{cor}[thm]{Corrolary}

\newtheorem{defn}{Definition}[section]
\newtheorem{ex}{Example}[section]
\newtheorem{rem}{Remark}[section] 
\newtheorem{con}{Condition}[section]
\endgroup

\newcommand\bthm {\begin{thm} \rm}
\newcommand\ethm { \end{thm}}
\newcommand\bprop {\begin{prop} \rm}
\newcommand\eprop { \end{prop}}
\newcommand\blem {\begin{lem} \rm}
\newcommand\elem { \end{lem}}
\newcommand\bcor {\begin{cor} \rm}
\newcommand\ecor { \end{cor}}
\newcommand\bdf {\begin{defn} \rm}
\newcommand\edf { \end{defn}}
\newcommand\bex {\begin{ex} \rm}
\newcommand\eex { \end{ex}}
\newcommand\brem {\begin{rem} \rm}
\newcommand\erem { \end{rem}}
\newcommand\bcon {\begin{con} \rm}
\newcommand\econ { \end{con}}

\newcommand\bitem { \begin{itemize}}
\newcommand\eitem { \end{itemize}}

\newcommand\beq {\begin{equation}}
\newcommand\eeq {\end{equation}}

\newcommand\bea{\begin{eqnarray}}
\newcommand\eea{\end{eqnarray}}

\newcommand\beaa{\begin{eqnarray*}}
\newcommand\eeaa{\end{eqnarray*}}

\newcommand{\Proof} {\noindent \it Proof . \rm}

\newcommand\bbib{\begin{thebibliography}{}}
\newcommand\ebib{\end{thebibliography}{}}

\newcommand\bref{\section*{References} 
\begin{list}{0}{\setlength{\rightmargin}{\leftmargin}}}

\newcommand\eref{\end{list}}

\newcommand{\sect}[1]{\section{#1}\setcounter{equation}{0}}

\newcommand{\muz}{\mu \quad \mbox{a.s.} \quad z}
\date{}
\title{Quantile regression in transformation models}
\author{Dorota M. Dabrowska \\
Department of Biostatistics \\
University of California \\
Los Angeles, CA 90095-1772}

\begin{document}
\setcounter{page}{1}

\maketitle

\begin{center}
Abstract
\end{center}

Conditional quantiles
provide a natural tool for reporting results from regression 
analyses  based on semiparametric transformation models.
We consider their estimation
and construction of confidence sets in the presence of censoring.

\it Key words: transformation models, censored data, quantile process. \rm

\sect{Introduction}
One-sided transformation models provide a popular tool for regression
analysis of failure time data. These models assume that the conditional
distribution of a failure time $T$ given a vector of covariates $Z$ has 
distribution function 
\beq
\bar F(t|z) = F(\Gamma(t), \theta|z) \; \quad \muz\;,
\eeq
where $\mu$ is the marginal distribution of covariates,
$\Gamma$ is an unknown increasing function mapping the support
of the marginal distribution of $T$ onto the positive half-line, and
$\F = \{F(x, \theta|z): \theta \in \Theta, x > 0\}$ is a parametric family of
conditional cdf's supported on $R^+$.
The most common choice  corresponds to the scale regression model
\beq
\bar F(t|z) = G(\Gamma(t) e^{\theta^T z}) \; \quad \muz,
\eeq
where $G$ is a known distribution function. In particular,
the proportional hazard model is of this form. In this case $G$ represents
exponential distribution and the unknown transformation $\Gamma$ is the 
so-called baseline cumulative hazard function. 
Proportionality of hazards means that the conditional
distribution of  $T$ given $Z=z$ has hazard rates $h(t|z)$
satisfying
$$
{e^{ -\theta^T z_1} \over e^{ -\theta^T z_2}} = 
{h(t|z_2) \over h(t|z_1)}
$$
for any two distinct covariate levels $z_1$ and $z_2$. 
This interpretation of parameters $(\Gamma, \theta)$ is lost in
other transformation models of type (2)
because the shape of the function $\Gamma$ depends on the 
distribution $G$. 

It is convenient to consider quantiles
$$
Q(p|z) = \inf \{t: \bar F(t|z) \geq p \}
$$
of the conditional distribution of $T$ given $Z=z$  as an alternative
parameter. 
In transformation models (2), we have
\beq
Q(p|z) = \Gamma^{-1}(e^{-\theta^T z} G^{-1}(p)) \quad 
\eeq
for all $p \in (0,1)$ and $\mu$ almost all $z$. Thus the conditional
quantiles are monotone in each coordinate of the vector 
$z = (z_1, \ldots,z_d)$. In addition,
the direction of monotonicity does not depend on $p$:
$$
\sign [{d \over dz_k} Q(p|z)] = \sign (-\theta_k)
\quad \for k = 1, \ldots,d \;.
$$
Invariance of the model with respect to the group of increasing
transformations implies also that 
for any $p_1 \not= p_2$ we have
\beq
{\Gamma(Q(p_1|z)) \over \Gamma(Q(p_2|z))} = 
{G^{-1}(p_1) \over G^{-1}(p_2)}  
\quad \muz
\eeq
and for any  $z_1 \not= z_2$
\beq
{\Gamma(Q(p|z_1)) \over \Gamma(Q(p|z_2))} = 
{e^{ -\theta^T z_1} \over e^{ -\theta^T z_2}} 
\eeq
for all $p \in (0,1)$. These three identities 
can be perhaps better understood by noting that (2) represents a
linear regression model 
$$
\log \Gamma(T) = -\theta^T Z + \epsilon \;,
$$
where $Z$ and $\epsilon$ are independent and $\exp \epsilon$ has distribution
function $G$. In linear regression models assuming that the transformation
$\Gamma$ is known and equal to $\Gamma(t) = t$, the conditional quantiles
 are linear
in $z$ but the slope of the regression does not change with $p$.
Likewise, the identities (4) and (5) have their additive 
analogue.  However, if the transformation is unknown, then 
the model is much more difficult to interpret in terms 
of the parameters $(\theta, \Gamma)$.   
 
Properties of quantile regression in the proportional
hazard model are further discussed in 
Koenker and Geling (2001) 
and Portnoy (2003). In particular, Koenker and Geling (2001)
proposed to measure the local effect of the 
regression coefficient on the conditional quantile $p$ in terms
of a  parameter $b(p, EZ) = [b_k(p, EZ), k = 1, \ldots,d]$, where
$$
b_k(p,z) = {d \over dz_k} Q(p|z) \;.
$$
This parameter can be applied to 
any regression model. In  (2) we have
$$
b(p, EZ) = - \theta^T {e^{-\theta^T EZ} G^{-1}(p) \over 
\gamma(Q(p|EZ))} \;,
$$
provided  the unknown transformation has
density $\gamma$ with respect to Lebesgue measure in a neighbourhood
of $Q(p|EZ)$.
While $b(p, EZ)$ is proportional to the regression coefficient 
$\theta$, the local
effect of the regression coefficient is determined 
by the shape of the density $\gamma$.   
Portnoy (2003)  considered direct modeling of the conditional
quantiles under the assumption that 
$\Gamma$  is the identity
map. His model takes form 
$$
 Q(p|z) = e^{\theta(p)^T z} \;, 
$$
so that for fixed $p$ the log-conditional quantiles are linear in $z$,
but also the quantile regression coefficient changes with $p$. 
 However, the
choice of the identity map may be problematic.
For other choices of the transformation, we have 
$Q(p|z) = \Gamma^{-1}(\exp{\theta(p)^T z})$.
Koenker and Geling's measure is given by
$$
b(p, EZ) = \theta(p)^T {e^{\theta(p)^T EZ} \over 
\gamma(Q(p|EZ))} \;.
$$
It shows that the model is more flexible than the
semiparametric transformation model (2), but it is 
not clear how to estimate the transformation function
in this setting.

In many practical situations researchers may be also interested
in the conditional distribution of $T$ given $\phi(Z)$, where 
$\phi$ is a known function. 
In particular, if $Z = (V,W)$ represents
a high-dimensional covariate, then the choice $\phi(Z) = V$ may
correspond to a low-dimensional vector of "main" covariates.
If $V$ and $W$ are dependent variables, then the conditional
distribution of $T$ given $V$ follows the more flexible 
transformation model (1).  For example, if (2) represents 
the proportional hazard model with parameters $\theta = (\theta_1, \theta_2)$
and the conditional distribution
of $\exp[\theta_2^T W] $ given $V$ is gamma with shape and scale 
equal to $\exp \xi(v)$
for a possibly nonlinear function $\xi$ of $v$, 
then the marginal conditional
distribution of $T$ given $V$ has distribution function of the 
form (1) with
$$
F(x, \theta_1, \xi |v) = 1 -
 \left (1 + \exp[\theta_{1} v + \xi(v)] x \right)^{-\exp[-\xi(v)]} \;.
$$
The ratio of conditional hazards is
$$
{h(x|v_2) \over 
h(x|v_1)} = {e^{ -\theta_1 v_1} \over   
e^{ -\theta_1 v_2}} \left[ { 1 + e^{\theta_1 v_1 + \xi(v_1)} \Gamma(x) \over   
1 + e^{\theta_1 v_2 + \xi(v_2)}  \Gamma(x)} \right ]
$$
For $x=0$ the right-hand side  is equal to $\exp [-\theta_1(v_1 - v_2)]$
and changes to  $\exp [\xi(v_1) - \xi(v_2)]$ as $x \uparrow \infty$.
 It represents an increasing function if $\theta_1(v_1 - v_2) \geq
 \xi(v_2) - \xi(v_1)$ and a  decreasing function,
 if the inequality is reversed.  
The conditional quantile function
is  equal to 
$$
Q(p|v) = \Gamma^{-1}(F^{-1}(p, \theta_1,\xi|v))
$$
where 
$$
F^{-1}(p, \theta, \xi|v) = 
\exp [-\xi(v)-\theta_1 v] [(1-p)^{-\exp \xi(v)} - 1] \;.
$$
If $\xi(v)$ is constant for almost all $v$, 
then we obtain the model (2). Otherwise the shape  of the 
quantile function changes with p. 
The ratios of the 
transformed quantiles (4) and (5) are no longer constant in $v$ and
$p$, respectively. 

In the general case, the conditional 
distribution of $e^{\theta_2^T W}$ given $V$ will not have a 
simple analytical form, even if specified via a parametric model.
However, 
quantile regression 
of the marginal conditional distributions of the 
failure time $T$ can also be estimated by combining nonparametric
regression with estimates of the parameters $(\theta, \Gamma)$.

In this paper we consider estimation of the conditional quantiles of 
 $T$ given $\phi(Z)$,
where $\phi$ is a function assuming a finite number of values.
In particular, if $Z = (Z_1, \ldots, Z_d)$ has one or more
discrete components, then results of this paper can be 
applied to estimation of quantiles of the marginal conditional distributions
of $T$ given any discrete  component of $Z$. On the other hand
in the case of continuous covariates estimation of the marginal 
conditional distribution and quantiles requires smoothing and may be
difficult to accomplish in  moderate  or heavily censored samples.
In such circumstances grouping observations into a small
number of categories provides an alternative.
For purposes of 
estimation of the parameters $(\theta, \Gamma)$ in transformation
models (1) and (2), we  use
procedures proposed by Bogdanovicius and Nikulin (1999) and Dabrowska (2005).
The approach allows for estimation of quantiles 
of the conditional distribution of $T$ given $Z=z$ much in the
same way as in the proportional hazard model, i.e. based
on the substitution  
of estimates of $(\theta, \Gamma)$ into (3) 
(Dabrowska and Doksum,
1987, Burr and Doss, 1993).
Here we  derive asymptotic structure of the 
estimates of the conditional quantiles under the assumption
that $\phi$ is a finite valued function, and  consider construction
of pointwise and simultaneous confidence sets. 
We also develop a
Gaussian  multiplier method for setting simultaneous confidence sets
for the conditional quantile function. It extends the
 Gaussian multiplier method
for setting
confidence bands for the conditional survival function in the proportional
hazard model (Lin, Fleming and Wei (1994)) 
to transformation models
 of type (1). In Section
3 we use data from a Vateran's Administration lung cancer clinical
trial (Kalbfleisch
and Prentice, 2000) to illustrate the results.
Section 4 contains proofs.

\sect{Estimation}
\setcounter{equation}{5}

We assume that the vector  $(X, \delta, Z)$
represents a nonnegative withdrawal time ($X$), a binary withdrawal 
indicator ($\delta = 1$ for failure and $\delta = 0$ for loss-to-follow-up)
and covariate ($Z$). The triple $(X, \delta, Z)$ is defined 
on a complete probability space $(\Omega, \F, P)$ and 
$(X, \delta)$ are given by 
$X = T \wedge \tilde T$, $\delta = 1(X = T)$, where $T$ and $\tilde T$
represent failure and censoring times. The variables $T$ and $\tilde T$
are 
conditionally independent given $Z$ and the conditional cumulative hazard
function of $T$ given $Z$ is of the form 
$$
H(t|z) = A(\Gamma_0(t), \theta_0|z) \quad \muz \;,
$$
where
 $\Gamma_0$
is an unbounded  continuous increasing function, 
 $\{A(x, \theta|z): \theta \in \Theta\}$ is a parametric 
family of cumulative hazard functions with hazard rate
$\alpha(u, \theta,z)$, and $\theta_0$ is the ``true'' parameter.
It is assumed throught the paper that the parameters
of the conditional distribution of the censoring times
are non-informative on $(\Gamma, \theta)$.

 Let
$N(t) = 1(X \leq t, \delta = 1)$ and $Y(t) = 1(X \geq t)$ denote
the counting and risk processes associated with the pair $(X, \delta)$.
We also set
$$
\tau_0 = \sup\{t: EY(t) > 0\}
$$
and assume the following regularity conditions.

\noindent
\bf Condition 1 \rm 

\bitem
\item[{(i)}] The covariate $Z$ has a nondegenerate distribution $\mu$
and is bounded: $\mu( |Z| \leq C) = 1$ for some constant $C$.
\item[{(ii)}] The function $EY(t)$ has at most a finite number of atoms, 
and $EN(t)$ is  continuous. 

\item[{(iii)}] The point $\tau > 0$ satisfies
$\inf \{t: E[N(t)|Z =z] > 0\} < \tau$ for $\muz$. In addition
$\tau < \tau_0$ if $\tau_0$ is a continuity point of the survival
function $EY(t)$, and $\tau = \tau_0$, if $\tau_0$ is an atom
of this survival function.

\item[{(iv)}]    
The parameter set $\Theta \subset R^d$ is open, 
and the   parameter $\theta$ is identifiable in the core model: $\theta \not= \theta' \tiff A(\cdot, \theta|z)
\not \equiv A(\cdot, \theta'|z) \quad \muz$.

\item[{(v)}] 
There exist constants $0 < m_1 < m_2 < \infty$ such that the hazard rate
$\alpha$ satisfies
\beq
m_1 \leq \alpha(x, \theta, z) \leq m_2
\eeq
for $\muz$ and all $\theta \in \Theta$, or 
(6) and (vi) holds for
$\tilde \alpha(x, \theta,z) = \alpha(\Phi(x), \theta,z) \Phi'(x)$, 
where $\Phi$
a strictly increasing unbounded twice continuously differentiable 
function $\Phi$ such that $\Phi(0) = 0$. 

\item[{(vi)}]
The function $\l(x, \theta,z) = \log \alpha(x, \theta,z)$ 
is twice continuously differentiable with respect to both $x$ and $\theta$.
The derivatives with respect to $x$ (denoted by primes) satisfy
$$
|\l'(x, \theta,z)| \leq \psi(x), \quad  |\l''(x, \theta,z)| \leq \psi(x) \;,
$$
where $\psi$ is a constant or a continuous bounded decreasing function.  
The derivatives with respect to $\theta$ (denoted by dots) satisfy
$$
|\dot \l(x, \theta,z)| \leq  \psi_1(x), \quad 
|\ddot \l(x, \theta,z)| \leq  \psi_2(x) 
$$
and
$$
|g(x, \theta,z) - g(x', \theta,z)| \leq \psi_3(x)[|x - x'| + 
|\theta - \theta'|] \;,
$$
where $g = \ddot \l, \dot \l'$ and $\l''$. The  functions 
$\psi_p, p =1,2, 3$ are continuous, bounded or strictly increasing 
and such that $\psi_p(0) < \infty$,
$$
\int_0^{\infty}  e^{-x} \psi_1^2(x) dx < \infty, \quad  
\int_0^{\infty}  e^{-x}  \psi_2(x) dx < \infty,  \quad
\int_0^{\infty}  e^{-x}  \psi_3(x) dx < \infty .
$$
\eitem

The assumption that the covariate  $Z$ is bounded is restrictive, but
standard for analysis of semiparametric models assuming that the 
transformation $\Gamma$ is unknown. In the special case
of the proportional hazard model, Andersen and Gill (1982) required
only existence of moments $EZ^2 e^{\theta^T Z} 1(X \geq x)$, for $x \geq 0$
in a neighbourhoood $\Theta \subset R^d$ of the true parameter $\theta_0$. 
However, setting $x = 0$, we see that this  moment condition
may lead to a constrained optimization problem which cannot be 
correctly stated, if the distribution $Z$ is unspecified. 
For example, if $Z$ is multivariate normal $N(0, \Sigma)$ and $\Sigma$ is
a known non-singular matrix, then 
the moment condition is satisfied for all $\theta \in R^d$ 
and the
usual unrestricted partial likelihood approach towards fitting the regression 
coefficients applies. 
However, if  $Z$ is a univariate lognormal variable, $Z \sim \exp \N(0,1)$,
then the parameter $\theta$ must be estimated under the 
added side condition $\theta \leq 0$. Thus the boundedness assumption
is restrictive, but allows for parameter estimation without
additional assumptions on the marginal distribution of the covariate.

Given an iid sample $(N_i, Y_i, Z_i), i = 1, \ldots, n$ of the 
$(N, Y,Z)$ processes, we set 
$N_.(t) = n^{-1} N_i(t)$,
$$
S(x, \theta ,t) = {1 \over n} \sum_{i=1}^n Y_{i}(t)
\alpha_{i}(x, \theta ) \;.
$$
and $\alpha_i(x, \theta)= \alpha(x, \theta|Z_i)$.
Following Bogdanovicius and  Nikulin (1999), define 
$$
\Gamma_{n\theta}(t) = 
\int_0^t {N_.(du) \over S(\Gamma_{n\theta}(u-), \theta,u)}\;, \quad 
\Gamma_{n\theta}(0-) = 0
$$
for any $\theta \in \Theta$.  
The process $\{\Gamma_{n\theta}: \theta \in \Theta\}$ 
is here thought as the sample analogue
of the Volterra integral equation
\beq
\Gamma_{\theta}(t) = \int_0^t {EN(du) \over s(\Gamma_{\theta}(u-),\theta,u)}
\;, \quad \Gamma_{\theta}(0-) = 0, \quad \theta \in \Theta \;,
\eeq
where $s(x, \theta,u) = EY_i(u) \alpha_i(x, \theta)$. 
The condition 1 (iv)  was used in Dabrowska (2005) to verify that this 
equation has a unique locally bounded solution, and such that 
$\Gamma_{\theta}(\tau_0) < \infty$ if $\tau_0$ is an
atom of the survival function $EY(t)$, and $\lim_{t \uparrow \tau_0}
\Gamma_{\theta}(t) \uparrow \infty$, if $\tau_0$ is a continuity point 
of $EY(t)$. In particular, the latter applies to   uncensored data.
Therein we show that in  the case of scale transformation models (2),
the condition 1 (v) is satisfied
by half-logistic, half-normal and half-t distributions, proportional
odds ratio distribution, frailty models with decreasing heterogeneity
with fixed frailty parameter  and polynomial hazards with nonnegative
constant coefficients. These models have smooth differentiable hazards
with respect to both $x$ and $\theta$ and integrability conditions
1 (vi) imply also that Fisher information is finite.
Affine independence of covariates 
is sufficient for the condition 1 (iv) to hold. In the case
of transformation models (1), the regularity conditions
are satisfied in the gamma frailty model
with frailty parameter representing a function of covariates
dependent on a Euclidean parameter. They are also satisfied
in  regular polynomial hazard regression models 
with nonnegative coefficients representing parametric  
functions of covariates.
In these models, the  conditional hazard rates 
are twice differentiable with respect to $x$,
while  the condition 1 (vi) imposes  a second order
differentiability assumption  on the functions of covariates.
Such differentiability conditions are in general not 
needed in regular parametric models. However, here we 
use semiparametric models and estimation of the parameter $\theta$
will be  based on a conditional rank statistics score equation. 
We do not know
at present time, how to relax these differentiability conditions 
to allow for estimation based on ranks.

 For any $\tau$ satisfying condition 1, 
the function $\{\Gamma_\theta(t): t \in [0, \tau],
\theta \in \Theta\}$ is Fr\'echet differentiable with respect to $\theta$
and the derivative satisfies the linear Volterra equation
$$
\dot \Gamma_{\theta}(t)  =   - \int_0^t \dot s(\Gamma_{\theta}(u-), \theta,u)
C_{\theta}(du) - \int_0^t \dot 
\Gamma_{\theta}(u-) s'(\Gamma_\theta(u-), \theta,u)
C_{\theta}(du) \;,
$$
where  $\dot s(\Gamma_{\theta}(u-),\theta,u) = EY_i(u) \dot 
\alpha_i(\Gamma_{\theta}(u-),\theta)$, 
$s'(\Gamma_{\theta}(u-),\theta,u) = EY_i(u) 
\alpha_i'(\Gamma_{\theta}(u-),\theta)$ and 
$$
C_{\theta}(t)  =  \int_0^t {EN(du) \over s^2(\Gamma_{\theta}(u-), \theta,u)}
\;. 
$$
In the case of the proportional hazard model,
the function $s'$ is identically equal to 0. 
Otherwise, the solution to this Volterra equation is given by
\beaa
\dot \Gamma_{\theta}(t) 
&  = &  - \int_0^t \dot s(\Gamma_{\theta}(u-), \theta,u) C_{\theta}(du) 
\P_{\theta}(u,t) \;, \\
\P_{\theta}(u,t) & = & \BIGP_{(u,t]} 
(1 - s'(\Gamma_{\theta}(w-),\theta,w) C_{\theta}(dw)) \;.
\eeaa
Here for any function $b$ of bounded variation, 
$\BIGP_{(u,t]} (1 + b(du))$ is the product integral, i.e.
$$
\BIGP_{(u,t]} (1 + b(dw))  =  \prod_{u < w \leq t} (1 + b(\Delta w)) 
\exp [b_c(t)]
$$
where $b_c$ is the continuous part of $b$ and the product is taken over
its atoms. 
To make the definition complete, in the 
case of the proportional hazard model we set $\P_{\theta}(u,t) \equiv 1$.
With this choice, the form of the function $\dot \Gamma_{\theta}$ is 
the same for all models of type (1) considered in this paper.

Let $\alpha_i(x, \theta) = \alpha(x, \theta, Z_i)$ and
$\l_i(x, \theta) = \log \alpha(x, \theta, Z_i)$. We shall apply
the same convention to derivatives of the functions
$\alpha_i$ and $\l_i$ with respect to $\theta$ and $x$.
Define functions
\beaa
\bar v(u, \theta) & = & { E Y_i(u) [\dot \l_i^{\otimes 2} \alpha_i]
(\Gamma_{\theta}(u), \theta) \over s(\Gamma_\theta(u), \theta, u)} -
\left ({\dot s \over s} \right)^{\otimes 2}(\Gamma_\theta(u), \theta, u) \\ 
v(u, \theta) & = & { E Y_i(u) [\l_i'^2 \alpha_i]
(\Gamma_{\theta}(u), \theta) \over s(\Gamma_\theta(u), \theta, u)} -
\left ({ s' \over s} \right )^2(\Gamma_\theta(u), \theta, u) \\
\rho(u, \theta) & = & 
{ E Y_i(u) [\dot \l_i  \l'_i \alpha_i]
(\Gamma_{\theta}(u), \theta) \over s(\Gamma_\theta(u), \theta, u)} -
\left ({\dot s \over s} \right )  \left ({ s' \over s } \right)
(\Gamma_\theta(u), \theta, u) 
\eeaa
and
\beaa
K_{\theta}(t,t') & = & \int_0^{ t \wedge t'} C_\theta(du) \P_\theta(u,t) 
\P_\theta(u,t') \\
B_{\theta}(t) & = & \int_0^t v(u, \theta) EN(du) \;.
\eeaa
Suppose that $v(u, \theta) \not \equiv 0$ a.e.--$EN$ and let
$\phi_{\theta} = \int_0^{\cdot} g_{\theta} d \Gamma_{\theta}$ be a vector
valued  function with $d$ components and 
square integrable with respect to $B_{\theta}$.

Define matrices  
\beaa
\Sigma_{1}(\theta) & = & \int_0^{\tau} 
 v_{\phi}(t, \theta) EN(du) \\
\Sigma_{2}(\theta) & = & 
 \int_0^{\tau} \int_0^{\tau} 
K_{\theta}(t,u) \rho_{\phi}(t, \theta) 
\rho_{\phi}(u, \theta)^T EN(du) EN(dt) \\
\Sigma(\theta) & = & \Sigma_1(\theta) + \Sigma_2(\theta)  
\eeaa
where
\beaa
v_{\phi}(t, \theta) & = & \bar v(t, \theta) + 
v(t, \theta) \phi_{\theta}^{\otimes 2}(t) 
- \rho(t, \theta) \phi_{\theta}^T(t)
 - \phi_{\theta}(t) \rho(t, \theta)^T \\
\rho_{\phi}(t, \theta) & = & \rho(t, \theta) - v(t, \theta) \phi_{\theta}(t)
\nonumber \;.
\eeaa
In the following we choose
 $\phi_{\theta}$ as solution to the Fredholm equation
\beq
\phi_\theta(t) + \int_0^{\tau} K_{\theta}(t,u) v(u, \theta) 
\phi_{\theta}(u) EN(du) = - \dot \Gamma_{\theta}(t) + 
\int_0^{\tau} K_{\theta}(t,u) \rho(u, \theta) EN(du) \;,
\eeq
or equivalently 
\beaa
\phi_\theta(t) + \dot \Gamma_{\theta}(t) & = &
\int_0^{\tau} K_{\theta}(t,u) \rho_{\phi}(u, \theta) 
EN(du) = \\
& = & \int_0^{\tau} K_{\theta}(t,u) \rho_{- \dot \Gamma}(u, \theta) 
EN(du) -
\int_0^{\tau} K_{\theta}(t,u) [\phi_{\theta}+\dot \Gamma_{\theta}](u)
B_{\theta}(du) \;.
\eeaa
This equation has a unique solution, square integrable with respect to
$B_{\theta}$. 
We define it as $\phi_{\theta} =
- \dot \Gamma_{\theta}$ if $\rho_{-\dot \Gamma}(u, \theta) \equiv 0$.
In this case we have $\Sigma_2(\theta) =  0$.
Finally, if $v(t,\theta) \equiv 0$ a.e. $EN$, 
then $\rho(t, \theta) \equiv 0$ as well.
For the sake of completeness we, set in this case $\phi_{\theta}= - \dot
\Gamma_{\theta}$. We also have $\Sigma_2(\theta)= 0$,
and $\Sigma_1(\theta)$ simplifies to $\Sigma_1(\theta) = \int \bar 
v(u,\theta) EN(du)$. This last choice corresponds to the proportional
hazard model, and  the scale regression models with regression
coefficient $\theta = 0$. (Note that if $v(u, \theta) \equiv 0$, then
the $\phi_{\theta}$ function does not enter into the score equation
below).

To estimate the parameter $\theta$, we use a solution to the 
score equation $U_n(\theta) = 0$, where
\beq
U_{n}(\theta) = {1 \over n} \sum_{i=1}^n \int_0^{\tau} 
[b_{1i}(\Gamma_{n\theta}(t), t, \theta)
- b_{2i}(\Gamma_{n\theta}(t), t, \theta)  \phi_{n\theta}(t) ] N_i(dt) \;,
\eeq
$\phi_{n\theta}$ is an estimator of $\phi_{\theta}$, and
$$
b_{1i}(x, t, \theta) =  \dot \l_i(x, \theta) - {\dot S(x, \theta,t) \over 
S(x, \theta,t)} \;,  \quad
b_{2i}(x, t, \theta)  =   \l_i'(x, \theta) - {S'(x, \theta,t) \over 
S(x, \theta,t)} \;.
$$
If $\Gamma_0$ is a known function, e.g. $\Gamma_0(t) = t$, then under
the assumption of conditional independence of failure and 
censoring times, the MLE score equation for
estimation of the parameter $\theta$ is given
by
$\tilde U_n(\theta) = 0$, where
$$
\tilde U_n(\theta) = {1 \over n} \sum_{i=1}^n
\int_0^{\tau} \dot \l_i(\Gamma_0(t), \theta) N_i(dt)
- \int_0^{\tau}  \dot S(\Gamma_0(t), \theta,t) \Gamma_0(dt)  
$$
and $\dot S(x, \theta,t) = n^{-1}\sum_{i=1}^n Y_i(t) \dot \alpha_i(x, \theta)$.
In addition, the assumption of conditional independence of failure
and censoring times implies that the function (7) satisfies
$\Gamma_{\theta_0}(t) = \Gamma_0(t)$ at the true value $\theta_0$ of 
the parameter $\theta$.  This last identity remains to hold
also when the transformation $\Gamma_0$ is unknown. Therefore a 
natural approach to estimation of the  parameter $\theta$ is to consider 
solving the score equation $\hat U_n(\theta) = 0$, where
$$
\hat U_{n}(\theta) = {1 \over n} \sum_{i=1}^n \int_0^{\tau} 
b_{1i}(\Gamma_{n\theta}(t), t, \theta) N_i(dt) \;.
$$ 
In particular, this is  the usual score equation
for estimation of the parameter $\theta$ in the proportional
hazard model.  In general transformation models (1), this 
choice leads to an asymptotically inefficient estimate
of the parameter $\theta$. It may also lead to estimates
of poor performance in moderate sample sizes. This also
applies to 
score processes of the form (9), where 
$\phi_{n\theta}$ is an estimate of some square integrable function
$\phi_{\theta}$ with respect to $B_{\theta}$. For example, Bogdanovicius 
and Nikulin (1999) considered the choice of $-\dot \Gamma_{\theta}$,
corresponding to the score equation derived 
from a modified partial likelihood function.
Under  mild regularity conditions on
the estimator of the the function $\phi_{\theta}$, the solution
to the score equation (9) exists with probability tending to 1
and is unique in local neighbourhoods of the true parameter $\theta_0$. 
However,  its asymptotic 
variance assumes the usual "sandwich" form because the 
process $\Gamma_{n\theta}$ has a non-trivial contribution 
to both asymptotic variance of the score process and the negative
derivative of it with respect to $\theta$. 
The choice of the $\phi_{\theta}$ function corresponding to 
the solution of to the Fredholm equation (8) leads to an $M$ estimator
whose asymptotic variance is of non-sandwich form and equal to
the inverse of the asymptotic variance of the score function.
The form of the solution to this equation can be found in Dabrowska
(2005).
The resulting estimator can also be shown to be asymptotically
efficient under the assumption that the point $\tau_0 =
\sup \{t: EY(t) > 0\}$ forms an atom of the survival
function $EY(t)$.
The following proposition summarizes some properties of the
estimates of $(\theta, \Gamma)$.

\bf Proposition 1 \rm Suppose that the conditions 1 are satisfied. 
 Let $\Sigma_1(\theta_0)$ be non-singular,
and let
 $\phi_{n\theta}$ be an estimator of this function
such that 
$\Vert  \phi_{n\theta_0} - \phi_{\theta_0} \Verti \to_P 0$,
$\limsup_n \Vert \phi_{n\theta_0} \Vertv = O_P(1)$,
$\phi_{n\theta} - \phi_{n\theta'} = (\theta - \theta') 
\psi_{n\theta, \theta'}$, where
$$
\sup\{ \limsup_n \Vert \psi_{n\theta,\theta'} \Vertv: 
 \theta \in  B(\theta_0, \epsilon_n) \} = O_P(1) 
$$
and
$B(\theta_0, \epsilon_n) = 
\{\theta: \Vert \theta - \theta_0 \Vert \leq \epsilon_n\}$ 
for some sequence  
$\epsilon_n \downarrow 0, \sqrt n \epsilon_n \to \infty$.
 Then, with probability tending to 1, the score equation
$U_{n}(\theta ) = 0$ has a unique solution 
 $\hat \theta$ in $B(\theta_0, \epsilon_n)$. 
Moreover,
$[\hat T, \hat W_0 ], \hat T = \sqrt n (\hat \theta - \theta_0)$,
$ \hat W_0 = \sqrt n [
\Gamma_{n\hat \theta} - \Gamma_{\theta_0} - (\hat \theta -
\theta_0) \dot \Gamma_{ \hat \theta}]$ converges weakly in 
$R^p \times \l^{\infty}([0,\tau])$ to a mean zero Gaussian process $[T, W_0]$
with covariance
\beaa
& & \cov T   =  \Sigma^{-1}(\theta_0)  \quad 
\cov (W_0(t),T)  = 
- \Sigma^{-1}(\theta_0) [\phi_{\theta_0}+ \dot \Gamma_{\theta_0}](t) 
\\ 
& & \cov(W_0(t), W_0(t'))   =   K_{\theta_0}(t,t') \;.
\eeaa

An example of an estimator of the function $\phi_{\theta}$ is 
given in Section 3. The asymptotic covariances can
be estimated using substitution method.

Let us assume now that $\D = \{D_j:j =1, \ldots,k\}$ is a finite partition
of the covariate space such that
\beq
\pi(D) = P(Z \in D) > 0, \quad D \in \D \;.
\eeq
We  denote by $F_D(t) = P(T \in t|Z \in D)$ the cdf of the conditional
distribution of  $T$ given $Z \in D, D \in \D$. Under the assumption of 
the transformation model, this function is of the form
$$
F_D(t) = {1 \over \pi(D)} E 1[Z \in D] F(\Gamma_0(t), \theta_0|Z) \;.
$$
In practice, the partition $D$ will be chosen based on the observations.
For example, if $Z= (Z_1, \ldots, Z_d)$ is a multivariate covariate, whose 
first component is continuous, then a natural partition of the 
covariate space may correspond to selection of $k=4$ intervals determined
by the sample quartiles of $Z_1$. If 
subjects are ranked according to values of 
the exponential factors $e^{\beta^T Z}$ than a natural 
partition may correspond to several groups 
determined by  the distribution of $e^{\beta^T Z}$. 
Any selection of such a partition requires some form of estimation
of parameters of the marginal distribution of the covariates. 
Here we consider a naive situation in which the cell probabilities 
can be estimated nonparametrically by means of sample proportions.
This choice arises in  analyses of models with possibly high-dimensional
discrete  or 
mixed discrete-continuous covariates, whenever 
interest is only in analyses of marginal conditional distributions 
corresponding to discrete variables  
representing treatment types,
patients' gender etc. In the data example given in section 3,
a many valued discrete variable representing a quantitative measurement 
patient's  performance status, admits  a natural partition 
into three groups corresponding to a more intuitive qualitative description 
of health condition at the time of entry into the clinical trial.

As an estimate  $\hat F_D(t)$ of the function $F_D(t)$ we take
\beaa
\hat F_D(t) & = & {1 \over \hat \pi(D)} {1 \over n} \sum_{i=1}^n 1(Z_i \in D)
F(\Gamma_{n\hat \theta}(t), \hat \theta|Z_i) \;,\\
\hat \pi(D) & = & {1 \over n} \sum_{i=1}^n 1(Z_i \in D) \;.
\eeaa
We also define scalar and vector valued functions
\beaa
\hat \psi_1(t,D) & = & {1 \over n} {1 \over \hat \pi(D)}\sum_{i=1}^n 1(Z_i \in D) 
f(\Gamma_{n\hat \theta}(t), \hat \theta|Z) \;,\\
\hat \psi_2(t,D) & = & \hat \psi_1(t,D) \dot \Gamma_{n \hat \theta}(t) +
{1 \over n} {1 \over \hat \pi(D)}
\sum_{i=1}^n 1(Z_i \in D) \dot F(\Gamma_{n\hat \theta}(t), 
\hat \theta|Z_i) \;,
\eeaa
where $\dot F(x,\theta|z)$ is the derivative 
of  $F(x,\theta|z)$ with respect to $\theta$. 

Finally, we denote by $\Vert \cdot \Vert$ the supremum norm on 
$\T = [0, \tau] \times \D$ and let $\l^{\infty}(\T)$ 
be the space of bounded functions on $\T$ endowed with the 
supremum norm. 

\bf Proposition 2 \rm
Suppose that the conditions of Proposition 1 are satisfied and (10) holds.

\bitem
\item[{(i)}] We have $\Vert \hat F - F \Vert \to_P 0$ and
$\hat W = \{\hat W(t,D) = \sqrt n [\hat F(t,D) - F(t,D)]: (t,D) \in \T 
\}$ converges weakly in $\l^{\infty}(\T)$ to $W$, a
mean zero Gaussian processes.
Its covariance function is given in Section 4.

\item[{(ii)}] Let $V_i = (V_{1i}, V_{2i}), i = 1,2, \ldots,n$ and $ V_{3} = 
(V_{31}, \ldots,
V_{3d})$ 
be mutually independent 
$\N(0,1)$ variables, independent of the observations $(X_i, \delta_i,Z_i),
i = 1, \ldots, n$. Define
\beaa
\hat W^{\#}_1(t,D) & = & {1 \over \sqrt n} {1 \over \hat \pi(D)}\sum_{i=1}^n V_{1n} 1(Z_i \in D)
[F(\hat \Gamma_{\hat \theta}(t), \hat \theta|Z_i) - \hat F_D(t)] \;,\\
\hat W_2^{\#}(t,D) & = & \hat W_0^{\#}(t) \hat \psi_1(t,D) + \int_0^{\tau}
\hat W_0^{\#}(s) \hat \rho_{\hat \phi_{n}}(s, \hat \theta) N_.(ds) 
\Sigma_n^{-1}(\hat \theta) \hat \psi_2(t,D) \;,\\
\hat W^{\#}_3(t,D) & = & V_3 \Sigma_{1n}^{1/2} (\hat \theta) 
\Sigma_n^{-1}(\hat \theta) \hat \psi_1(t,D) \;, 
\eeaa
where
$$
\hat W_0^{\#}(t) = {1 \over \sqrt n} \sum_{i=1}^n V_{2i} {1[X_i \leq t, 
\delta_i = 1] \over S(\Gamma_{n \hat \theta}(X_i-), \hat \theta, X_i)} 
\P_{n\hat \theta}(X_i, t)
$$
and $\Sigma_{1n}(\hat \theta)$,  $\Sigma_{n}(\hat \theta)$, $\hat \phi_{n} =
\phi_{n\hat \theta}$, $\hat \rho_{\hat \phi_n}(u, \hat \theta)$ and, 
$\P_{n\hat \theta}(u,t)$ are estimates of $\Sigma_1(\theta_0)$, 
$\Sigma(\theta_0)$,
$\phi_{\theta_0}$, $\rho_{\phi_{\theta_0}}(u, \theta_0)$, $\P_{\theta_0}(u,t)$
obtained using substitution method. The  process
$\hat W^{\#} = \{\hat W^{\#}(t,D) = \sum_{j=1}^3
 \hat W^{\#}_j(t,D):
(t,D) \in \T\}$ converges weakly (unconditionally) in $\l^{\infty}(\T)$ to a
Gaussian process $W^{\#}$  with
 the same covariance function as the process
$W$ of part (i) and independent of it. Conditionally, 
the process $\tilde W$ converges weakly to $W$ in probability.  
\eitem

The proof is given in Section 4.
In the first part of the proposition, the  
observations
$R_i = (X_i, \delta_i,Z_i),
i = 1, \ldots, n, \ldots$. 
are defined
as coordinate projections on the 
product probability space $(\Omega^{\infty}, \F^{\infty},
P^{\infty})$. In the second part, 
we use the product probability space
$(\Omega^{\infty} \times \V \times \V', \F^{\infty} \times
\B \times \B', P^{\infty} \times Q \times Q')$. The variables  
$R_i = (X_i, \delta_i,Z_i),
i = 1, \ldots, n, \ldots$, $V_i, i = 1, \ldots,n \ldots$ and $V_3$
are defined as first, second and last projections.
Conditional weak convergence in probability
means 
$$
\sup_{f \in BL_1} |E_V^* f(W^{\#}) - Ef(W)| \to 0
$$
in (outer) probability, where $BL_1$ is the set of all real functions
on $\l^{\infty}(\T)$ with a Lipschitz norm bounded by 1  
(van der Vaart and Wellner, 1996, Ch. 2.9).

We proceed to the discussion of the properties of the quantile regression.
For $p \in (0,1)$ and (fixed ) $D \in \D$ let
$$
\l_D(p) = \inf\{t: F_D(t) \geq p\}\;, 
\quad u_D(p) = \sup\{t: F_D(t) \leq p\} \;.
$$
Then $\l_D(p) \leq u_D(p)$  and the $p$-th quantiles of the conditional
distribution of $T$ given $Z \in D$ are defined as the set of numbers
in the closed interval $[\l_D(p), u_D(p)]$. We denote by 
$\hat \l_D(p)$ and $\hat u_D(p)$ the sample counterparts of  these 
points, i.e. 
$$
\hat \l_D(p) = \inf\{t: F_D(t) \geq p\}, \quad 
\hat u_D(p) = \sup\{t: F_D(t) \leq p\} \;.
$$
 If 
$u_D(p) < \tau$, then under assumptions of Proposition
2, we have 
\beq
\l_D(p) \leq \liminf_n \hat \l_D(p) \leq \limsup_n \hat u_D(p) \leq u_D(p) 
\eeq
with probability tending to 1. Indeed, let $\epsilon = \epsilon(D) > 0$ 
be arbitrary but 
small enough so that $u_D(p) + \epsilon < \tau$. Then
$$
F_D(\l_D(p) - \epsilon) < p, \quad  F_D(u_D(p) + \epsilon) > p 
$$ 
and uniform consistency of the estimate $\hat F_D(\cdot)$ implies that
with probability tending to 1, we also have
$$
\hat F_D(\hat \l_D(p) - \epsilon) \leq p, \quad  
\hat F_D(\hat u_D(p) + \epsilon) \geq  p \;. 
$$
This in turn implies (11).

In the following we shall assume that the transformation function $\Gamma_0$
has density $\gamma$ with respect to the Lebesgue measure, and the function
$\gamma$ is uniformly continuous and bounded away from 0 on an
interval $[0, \tau_1 - \epsilon, \tau_2+\epsilon], 0 < \tau_1 - \epsilon, \tau_2 + \epsilon \leq \tau \leq \tau_0$ and such that 
\beq
\tau_1 = \min \{ \l_D(p_1): D \in \D\}, \quad \tau_2 = 
\max \{ u_D(p_2): D \in \D \} \;.
\eeq
Let $I = [p_1, p_2]$ and  set  $\I = I \times \D$. 
In this case the conditional distribution of $T$ given $Z \in D$ has a unique
$p$-th quantile $Q_D(p)$ for any $p \in I$ and we define its sample
analogue by setting
$$
\hat Q_D(p) = \hat \l_D(p) = \inf\{t: \hat F_D(t) \geq p\} \;. 
$$
Then (11) implies that $\hat Q_D(p) \to_P Q_D(p)$ pointwise in 
$(p, D) \in \I$. Using finiteness of the class $\D$, monotonicity of 
$F_D(t)$ and $\hat F_D(t)$, and an argument similar to the classical
Glivenko-Cantelli theorem, 
we also have 
$$
\sup \{| \hat Q_D(p) - Q_D(p)|: (p,D) \in \I\} \to_P 0 \;.
$$

\bf Proposition 3 \rm
Suppose that the conditions of Proposition 2 hold, and
$\Gamma_0$
has density $\gamma$ with respect to the Lebesgue measure such that 
$\gamma$ is uniformly continuous and bounded away from 0 on an
interval $[0, \tau_1 - \epsilon, \tau_2+\epsilon], 0 < \tau_1 - \epsilon, \tau_2 + \epsilon \leq \tau$ satisfying  (12). The normalized quantile
process  $\hat V = \{\hat V(p,D): (p,D) \in \I \}$ given by
$$
\hat V(p,D) = \sqrt n[ \hat Q_D- Q_D](p) \;, 
$$
converges weakly in $\l^{\infty}(\I)$ to 
  $V = \{ V(p,D) = - h(p,D) W(Q_D(p),C): (p,D)  \in \I\}$, 
where
$$
h(p,D) = [f_D(Q_D(p)) \gamma(Q_D(p))]^{-1} \;.
$$

\Proof 
We have
$\hat V(p, D) = \hat h(p,D) \hat R(p,D)$,
where 
\beaa
\hat h(p,D) & = & \left ({\hat Q_D - Q_D \over F_D \circ \hat Q_D - F_D \circ Q_D}
\right ) (p) \;,\\
\hat R(p, C) & = & \sqrt n [F_D \circ \hat Q_D - F_D \circ Q_D](p)
\;.
\eeaa
Since the function $\gamma$ is positive and uniformly continuous on
$[\tau_1 - \epsilon, \tau_2 + \epsilon]$, uniform consistency of the 
sample quantile function implies
$$
\sup \{|\hat h - h|(p,D): (p,D) \in \I \} \to_P 0 \;.
$$
The process $\hat R(p,D)$ is on the other hand given by
$\hat R(p,D) = \sum_{j=1}^3 \hat R_j(p,D)$, where 
\beaa
\hat R_1(p, D) & = & - (\hat W_D \circ Q_D)(p) \;, \\
\hat R_2(p, D) & = & - (\hat W_D \circ \hat Q_D - \hat W_D \circ Q_D)(p) \;,\\
\hat R_3(p,D) & = & \sqrt n [\hat F_D \circ Q_D(p) - p] \;.
\eeaa
We have $\sup\{|\hat R(p,D)|: (p,D) \in \I \} \leq 
\sup \{ |\hat W_D(u) - \hat W_D(u-)|: u \in
[\tau_1 - \epsilon, \tau_2 + \epsilon], D \in \D\} = 
O_p(n^{-1/2})$ because the 
function $\hat F_D(x)$ has jumps of order $O_p(n^{-1})$. Application
of the Skorohod-Dudley-Wichura construction implies also that
$\sup \{| \hat R_2(p, D)|: (p,D) \in \I\} \to_P 0$, while the 
process  $\{ \hat R_1(p, D): (p,D) \in \I\}$ 
converges weakly in $\l^{\infty}(\I)$
to $\{- W_D \circ Q_D(p): (p,D) \in \I \}$.  \done

We shall apply now this result to construct pointwise confidence 
intervals for the $p$-th quantile.  Let $v_D(t)$ be the asymptotic
variance function
of the process $\{W(t,D): (t,D) \in \T\}$. It is derived in Section 4.
Here we shall use only that this function is positive and continuous
on the interval $[\tau_1 - \epsilon, \tau_2 + \epsilon]$, and its 
its plug-in analogue $\hat v_D(t)$ is uniformly consistent on the set
$[\tau_1 - \epsilon, \tau_2 + \epsilon] \times \D$.

For $p \in (0,1)$ and $D \in \D$, let
$$
p_{n}^{\pm} = p \pm {1 \over \sqrt n} \hat v_D(\hat Q_D(p)) z(\alpha) \;,
$$
where $z(\alpha)$ is the upper $\alpha/2$ percentile of $\N(0,1)$ distribution.
Proposition 3 and the inequalities
\beaa
& & \hat Q_D(p) \geq s \tiff p \geq \hat F_D(s) \;,\\
& &  Q_D(p) \geq s \tiff p \geq  F_D(s) \;,
\eeaa
imply that $[\hat Q_D(p_n^-), \hat Q_D(p_n^+)]$ is a
$100\% \times (1-\alpha)$ asymptotic pointwise confidence interval for the 
conditional quantile $Q_D(p)$. 

Unfortunately, in practice the points
$p_n^{\pm}$ may fall outside the range $[0,1]$. To circumvent this
problem, we follow the approach of Bie \it et al. \rm (1987) and 
consider confidence intervals based on transformations. Let
 $g$ be a strictly monotone cdf with density $g'$ supported on the 
whole real line. Set
$$
p_{nD}^{\pm} = 
g^{-1}(p) \pm {1 \over \sqrt n} {\hat v_D(\hat Q_D(p)) \over 
g'(g^{-1}(p))} z(\alpha) \;.
$$
With probability tending to 1, the inequalities 
$$
\hat Q_D(g(p_{nD}^{-})) \leq Q_D(p) \leq \hat Q_D(g(p_{nD}^+))
$$
are equivalent
to
$$ 
-z(\alpha) \leq g'(g^{-1}(p)) \sqrt n {g^{-1}(\hat F_D (Q_D(p))) 
- g^{-1}(p) \over \hat v_D(\hat Q_D(p))} \leq z(\alpha) 
$$
and application of delta method implies that
$[\hat Q_D(p_{nD}^-), \hat Q_D(p_{nD}^+)]$ is a
$100\% \times (1-\alpha)$ asymptotic confidence interval for the 
conditional quantile $Q_D(p)$.

Construction of simultaneous confidence sets for the function
$\{Q_D(p): (p,D) \in \I\}$ is more difficult because the process
W appearing in Propositions 2 and 3 forms a sum of independent Gaussian
processes with correlated increments. Therefore, 
following Burr and Doss (1993) and
Lin, Fleming and Wei (1994), we propose the use of simulated confidence sets.

Define 
$$
U = \sup \{{|W(Q_D(p),C)| \over v_D(Q_D(p))}: (p,D) \in \I \} =
\sup \{{|W(t,D)| \over v_D(t)}: t \in [Q_D(p_1), Q_D(p_2)], D \in \D \} 
$$
and let $u(\alpha)$ be the upper $100\%(1-\alpha)$
percentile of its distribution. To obtain an approximation to the 
critical level $u(\alpha)$, we generate  mutually independent
standard normal vectors $V$ defined as in Proposition 3, and form 
$$
U^{\#} = \sup \{{|\hat W^{\#}(t,D)| \over \hat v_D(t)}:
t \in [\hat Q_D(p_1), \hat Q_D(p_2)], D \in \D\}
$$
The procedure is repeated independently $m$ times, for some large $m$,
to obtain $m$ iid copies $U_1^{\#}, \ldots, U_m^{\#}$. The estimate $u^{\#}(\alpha)$
of the critical point $u(\alpha)$ is taken as the empirical
$(1-\alpha)$ quantile of $U_1^{\#}, \ldots, U_m^{\#}$. The corresponding
simulated confidence set for $\{Q_D(p): (p,D) \in \I\}$ is chosen as
$$
\{[\hat Q_D(\hat p_{nD}^-), \hat Q_D(\hat p_{nD}^+)]: D \in \D\} \;,
$$
where
$$
\hat p_{nD}^{\pm} = 
g^{-1}(p) \pm {1 \over \sqrt n} {\hat v_D(\hat Q_D(p)) \over 
g'(g^{-1}(p))} u^{\#}(\alpha) \;.
$$
Application of Propositions 2-3
implies that $u^{\#}(\alpha)$, the  upper $\alpha$--quantile of this (conditional)
distribution satisfies $u^{\#}(\alpha) \to u(\alpha)$ in probability.

An alternative approach to construction of simultaneous confidence
sets may be based on bootstrap. Lin, Fleming and Wei (1994) argued
that in the case of Cox regression with external time dependent
covariates, it is not clear how to implement bootstrap to
construct simultaneous confidence bands for the conditional survival
function, or other functionals related to it. In our setting
covariates are time independent, and confidence sets can be based
on ``obvious'' bootstrap. We can draw $R^{*}_n = [(X_i^{*},
\delta_i^{*}, Z_i^{*}): i = 1, \ldots,n]$  by sampling with
replacement from the empirical distribution function of the 
$[(X_i,
\delta_i, Z_i): i = 1, \ldots,n]$ observations For each sequence
$R_{nj}^*: j = 1, \ldots, m$  
we can compute  bootstrap estimates $\{Q^{*}_D(p), (p,D) \in \I\}$
and next use them to approximate the distribution of the 
quantile process. Although it is possible to show consistency
of this procedure, its
drawback  lies in the 
computational burden needed to
construct estimates  $(\theta^{*}_n, \Gamma_{n\theta^{*}}^{*})$ 
for each of the $m$
simulated data sets. 
In the case of the proportional hazard model, Hjort (1985) proposed
the use of ``model based'' bootstrap. Burr and Doss (1993)
applied it to the construction of simultaneous confidence bands 
for the conditional
median. In this approach, the distribution of the quantile process
is approximated based on  artificial observations 
$(X^{*}_i, \delta^{*}_i), i = 1, \ldots, n$ defined as 
$X^{*}_i = T^{*}_i \wedge \tilde T^{*}_i, \delta^{*}_i = 1(T^{*}_i \leq \tilde T^{*}_i)$,
where $T_i^{*}$ is sampled from the distribution 
$F(\hat \Gamma_{n\hat \theta}(t), \hat \theta|Z_i)$ and $\tilde T_i^*$
is sampled from $\hat G(t) = 1 -$ Kaplan-Meier estimate of the censoring
distribution.  This approach uses the assumption
that censoring time is independent of covariates, which need not be
satisfied in many practical situations. It is in principle
possible to relax it  by chosing a parametric or a 
semi-parametric
model for the conditional distribution of censoring times, however, selection
of such a model is often quite difficult, and its misspecifaction  
may affect
the performance of confidence procedures.

\sect{Example}
\setcounter{equation}{12}

For illustrative purposes we consider now data from the Veteran's 
Administration lung cancer trial (Kalbfleisch and Prentice, 2000). 
In this trial
males with inoperative lung cancer were randomized to either a 
standard or an experimental chemotherapy treatment and subsequently
followed until death or withdrawal from the study. 
We shall look at the subgroup 
of 97 patients, who received no prior therapy, and use
two covariates
corresponding to 
performance   status at the time of entry into the clinical trial 
and histopathological type of tumor (squamous, small cell, adeno and 
large cell). 

Several authors (e.g. 
Bennett ( 1983), Pettit (1984), Cheng et al.  (1995) and Murphy, Rossini
and van der Vaart  (1996)) 
 proposed the use of the proportional odds ratio for analysis of this 
dataset. Our estimates are easy to compute in this case 
because the hazard rate of the $i$-th subject satisfies
\beq
\alpha_i(x, \theta) = e^{\theta^T Z_i}(1 + e^{\theta^T Z_i} x)^{-1}, 
\quad \l_i'(x, \theta) = - \alpha_i(x, \theta), \quad \dot \l_i(x, \theta) =
Z_i e^{-\theta^T Z_i} \alpha_i(x, \theta) \;.
\eeq
For fixed $\theta$, the estimate $\Gamma_{n\theta}$ is computed based
on the recurrent formula given by Bogdanovicius and Nikulin (1999):
$$
\Gamma_{n\theta}(t) = \Gamma_{n\theta}(t-) + {N_.(\Delta t) \over
S(\Gamma_{n\theta}(t-), \theta,t)}
$$
with the initial condition  $\Gamma_{n\theta}(0-) = 0$. 
The sample version of the function $\dot \Gamma_{\theta}$ can
be evaluated as
$$
\dot \Gamma_{n\theta}(t) = 
\dot \Gamma_{n\theta}(t-) -
[\dot S(\Gamma_{n\theta}(t-), \theta,t) +      S'(\Gamma_{n\theta}(t-) 
 \theta,t) \dot \Gamma_{n\theta}(t-)] {N_.(\Delta t) \over
S^2(
\Gamma_{n\theta}(t-), \theta,t)}
$$
and  $\dot \Gamma_{n\theta}(0-) = 0$. 
The solution to the Fredholm equation can be obtained
as follows. Let
$X_{(1)} < \ldots < X_{(m)}, m \leq n$ be the 
distinct uncensored observations in the sample. 
Dropping dependence on the parameter $\theta$, let
$B_n, C_n$ be the plug-in sample analogues of the functions 
$B_{\theta}$ and
$C_{\theta}$. These are step functions
with jumps at points $X_{(i)}$ and we arrange their jumps
into $m \times m$ diagonal matrices 
 ${\bf B_n(\Delta X)} = \diag \{B_n(\Delta X_{(i)}): i = 1, \ldots, m \}$,
and ${\bf C_n(\Delta X)} = \diag \{C_n(\Delta X_{(i)}): 
i = 1, \ldots, m \}$. let
 $\bf \rho_n(X)$ be an $m \times d$ matrix of the sample 
analogues of the conditional covariances $\rho_{-\dot \Gamma}(u, \theta)$
at points $X_{(i)}, i = 1, \ldots, m$.
(Here $d$ is dimension  of the parameter $\theta$). The matrix 
${\bf C_n(\Delta X)}$ has positive entries, the matrix
${\bf B_n(\Delta X)}$ nonnegative. 
If ${\bf B_n(\Delta X)} 
 \equiv 0$ then also ${\bf \rho_n(X)} \equiv 0$.
Setting $\psi_{n\theta} = \phi_{n\theta} + \dot \Gamma_{n\theta}$,
the 
discrete version of the Fredholm equation
corresponds to 
$$
\bf{[I +  K_n(X) B_n(\Delta X)] 
\psi_n(X)} = \bf{K_n(X)\rho_n(X)} \;,
$$
where  ${\bf \psi_n(X)} = 
[\psi_n(X_{(i)}): i = 1, \ldots,m]^T$ is an $m \times d$ matrix of
unknowns,  
$\bf K_n(X)$
is an $m \times m$ matrix with 
entries ${\bf K_n(X)} = [K_n(X_{(i)}, X_{(j)})]$ and $\bf I$
represents an $m \times m$  identity. If $\bf B_n(\Delta X) \equiv 0$ or 
$\bf \rho_n(X) \equiv 0$ then the solution is $\bf \psi_n(X) \equiv 0$.
Otherwise, 
${\bf \psi_n(X)  =  
P_n^T(X) g_n^{-1}(X) 
P_n(X) \rho_n(X)}$, 
where
${\bf g_n(X)} = [g_{ij}]$ 
is a tridiagonal symmetric matrix with entries
$g_{ii} = c_{i} + c_{i+1} + b_i, g_{i, i+1} = - c_{i+1} = g_{i+1,i}, i =1, \ldots, m-1$ and $g_{mm} = c_m + b_m$, where $b_i = \P_{n\theta}(0, X_{(i)})^2 
B_n(\Delta X_{(i)}),  b_i = \P_{n\theta}(0, X_{(i)})^2 
C_n(\Delta X_{(i)})^{-1}$ and
${\bf P_n(X)} = 
  \diag [\exp ( -\int_{[0, X_{(i)}]}  
S'(\Gamma_{n\theta}(u-), \theta, u) C_{n\theta}(du)): i 
=1, \ldots, m] \sim 
\diag [ \P_{n\theta}(0, X_{(i)}), i = 1, \ldots,m]$.
(Dabrowska, 2005). After obtaining the solution,
$\psi_{n\theta}$ we set $\phi_{n\theta} = \psi_{n\theta} - 
\dot \Gamma_{n\theta}$. 
The estimate $\hat \theta$ can be obtained  using Fisher
scoring algorithm. The
algorithm can be started by setting $\hat \theta^{(0)}$ obtained
by solving the same score equation, but function $\phi_{n\theta}$ set to 0
or $-\dot \Gamma_{n\theta}$.

The estimate $ \Gamma_{n\hat \theta}$ is a cadlag step function
with jumps at uncensored observations, and so is the estimate $\hat
F_D(t)$ of the conditional distribution 
function of $T$ given $Z \in D$. Thus the graph of the quantile 
function can be obtained by inverting graphically the plot of this function.
The estimate $\hat v_D(t)$ of the asymptotic variance of the
$\sqrt n [\hat F_D -
F_D](t)$ and  the process $\hat W^{\#}(t,D)$ can be easily computed
based on expressions  given in Sections 2 and 4.

Table 1 provides regression coefficients and their standard errors
for the Veteran's Administration lung cancer data. 
In this data set the performance score (PS) 
has range between 10 and 99, with lower values
indicating poorer performance status at the time of entry
into the trial. This covariate was used in the regression
model after standardizing it to have average zero and standard deviation 1. 
The negative sign of the  regression coefficient 
indicates that patients
with higher performance score have lower odds on death and thereby a better 
survival experience. Patients with squamous tumor have a slightly
lower odds on death than large cell tumor patients, however,
the difference  
is not significant. Patients with  adeno or small cell tumor have
higher odds on death than patients with squamous or large cell types.

\begin{center}
Table 1 about here
\end{center}

We  shall consider now two partitions $\D$ of the covariate space.
In both cases, we shall consider quantile regression
estimates in the range $p \in (.25, .75)$. 
Simultaneous confidence sets
are  based on the transformation 
$g^{-1}(p) = 
\log( -\log (1-p))$ and 
we used  1000 Monte Carlo simulations of the V vectors 
(section 2) to obtain the critical points.

The first partition corresponds to the four histopathological types of
tumor. Figure 1 shows the corresponding
quantile regression and confidence set
for the  conditional quantiles.
The plots support results of Table 1 and show that patients
with  squamous or large cell tumor  perform better than patients with 
adeno or small tumor cells. However, within each pair of tumor types,
the confidence sets are nearly the same so that the differences are small.

\begin{center}
Figure 1 about here
\end{center}

Next we partition the covariate space according to the performance 
status at the time of entry into the trial. We consider patients, who
are  completely 
hospitalized  (PS $< 40$), partially confined  
(PS $\in [40,70)$) and  who are not able to care 
(PS $\geq 70$). In 
Figure 2, the confidence sets for the hospitalized and partially
confined patients nearly overlap, suggesting similar survival experience
after treatment. This experience is much worse than for patients who 
are not able to care. For example, the estimated median time till death
for hospitalized, partially confined and unable to care patients
is  25, 29 and 110 days, respectively. The corresponding confidence
bounds are $(22, 35)$, $(24, 36)$ and $(103, 112)$ days.
Figure 3.2 suggests also that effect of the PS score is not linear,
and a regression model using a binary covariate: $Z = 1 (0)$ if PS score 
$\geq   (<)   70$ may be more appropriate.

\begin{center}
Figure 2 about here
\end{center}

We have also considered the choice of the proportional hazard
model and generalized inverse Gaussian frailty model. 
In each of these models the regression coefficients had the same
sign, however, neither of the transformation models could be 
fully justified. In Figure 3 we show  nonparametric plots of 
the Aalen-Nelson
estimator, odds ratio function and Kaplan-Meier estimator of the 
survival function for the four tumor cell types : squamous c
(solid line), large  (dotted line), small  (short dash)
and adeno  (log dash). The plots of 
the cumulative hazard function of the large and squamous cell
type cross at around 150 days. Patients with squamous 
cell type are initially at a higher risk for death but at around 150
days after treatment the role of the two groups is reversed. The 
corresponding plots of the odds ratio function suggest that 
the choice of proportional hazard model may not be appropropriate
and that odds ratio functions are close for the two groups. 
In the case of the adeno and small cell tumor cell type groups, the 
graphs of both cumulative hazard and odds ratio
functions cross only at the upper tail, however, the two groups
can be only compared during the initial 180 days. 

\begin{center}
Figure 3 about here
\end{center}

These graphs illustrate typical difficulty arising
in regression analyses based on transformation models of type (1)
or (2).
The transformation models  assume that the  conditional 
distributions of the failure time $T$ given $Z =z$
 have the same support as the marginal distribution of $T$
for $\mu$-  almost all $z$. 
This assumption fails to be satisfied in the fully
nonparametric setting, not assuming any restrictions on the 
support or shape of the conditional distribution of $T$ given $Z =z$.
If $\bar F(t|z)$ represents the conditional distribution
function of $T$ given $Z =z$ and $\bar G$ is the corresponding
marginal distribution function of $T$, then 
setting
\beaa
\tau_1(z) = \inf\{t: \bar F(t|z) > 0 \} & & 
\tau_2(z) = \sup\{t: \bar F(t|z) < 1 \} \\
 \tau_1 =  
\inf\{t: \bar G(t) > 0 \}  & & 
\tau_2 = \sup\{t: \bar G(t) < 1 \} 
\eeaa
we have $\tau_1 \leq \tau_1(z) \leq \tau_2(z) \leq \tau_2$ for 
$\mu$-almost all $z$, i.e. the marginal distribution
of $T$ has longer support than the conditional distributions. 
For different covariate levels $z_1$ and $z_2$,
the intervals $[\tau_1(z_1), \tau_2(z_1)]$ and
$[\tau_1(z_2), \tau_2(z_2)]$ may be very different. 

In the present
example, 
large and squamous cell type patient groups have longer support interval
than the groups of squamous and adeno cell types. Apparently,
patients for whom treatment is beneficial live longer. The choice
of the proportional odds ratio model  
appears to be more  appropriate than the 
proportional hazards model, however, 
it does not accommodate 
variable support intervals of conditional distributions of different
subgroups. The problem applies to all transformation models
of type (1) and (2). The plots of Kaplan-Meier
estimators corresponding to the four groups are proper
survival functions in this data example because data are lightly
censored (Kalbfleisch and Prentice, 2000). 
In moderately or heavily censored samples, the grouped data
Kaplan-Meier estimator will often form an improper survival
function. In such circumstances,  variable 
supports of Kaplan-Meier estimator may indicate also presence
of informative censoring. The difficuties in handling
variable supports of conditional distributions apply
also to other common parametric and semiparametric
regression models in survival analysis and are very common
in practical applications.

\sect{Proofs}
\setcounter{equation}{13}

In this section, we denote by $M_i(t)$ the process
$$
M_i(t) = 1(X_i \leq t) - \int_0^t Y_i(u) \alpha_i(\Gamma_{\theta_0}(u), \theta_0)
\Gamma_{\theta_0}(du) \;,
$$
where $\Gamma_0 = \Gamma_{\theta_0}$ is the ``true'' transformation.
Then $M_i$ are independent mean zero martingales, with respect to natural
filtration generated by $\F_t = \sigma\{(N_i(s), Y_i(s+), Z_i): s \leq t, 
i = 1, \ldots, n\}$.
 For any measurable functions $g_q(u,z), q = 1,2$ such that
$$
E\int Y_i(u)g^2_q(u,Z_i) \alpha_i(\Gamma_0(u),\theta_0)\Gamma_0(du) < \infty
$$ 
we have
\beaa
& & \cov (\int g_1(u,Z_i) M_i(du), \int g_2(u,Z_i) M_i(du)) = \\
& &
E\int Y_i(u)g_1(u,Z_i) g_2(u,Z_i) \alpha_i(\Gamma_0(u),\theta_0)\Gamma_0(du) 
\eeaa

\bf Lemma 1 \rm
 Suppose that the conditions of Propositions 1 and 2 are 
satisfied.
\bitem
\item[{(i)}] 
The estimate $\hat \theta$ satisfies $\sqrt n[\hat \theta - \theta_0] =
\Sigma(\theta_0)^{-1} \sqrt n \tilde U_n(\theta_0) + o_P(1)$, where
$\Sigma(\theta) = \Sigma_1(\theta) + \Sigma_2(\theta)$ and 
$\tilde U_n(\theta_0) = n^{-1}\sum_{i=1}^n [U_{1i}(\theta_0) + 
U_{2i}(\theta_0)]$
is given by
\beaa
U_{1i}(\theta_0) 
& = & \int_0^{\tau} b_i(\Gamma_{\theta_0}(u), \theta_0,u) M_i(du) \;,\\
U_{2i}(\theta_0) & = & - \int_0^{\tau} W_{0i}(t)
    \rho_{\phi_{\theta_0}}(t, \theta_0) EN(dt) \;,\\
W_{0i}(t)  & = & \int_0^t {M_i(du) \over s(\Gamma_0(u-), \theta_0,u)}
 \P_{\theta_0}(u,t) 
\eeaa
and
\beaa
b_i(\Gamma_{\theta_0}, \theta_0,u) & = & 
\dot \l_i(\Gamma_{\theta_0}(u), \theta_0) -
\l'_i(\Gamma_{\theta_0}(u), \theta_0) 
\phi_{\theta}(u) \\
& - & {\dot s \over s}(\Gamma_{\theta_0}(u), \theta_0,u) +
{ s' \over s}(\Gamma_{\theta_0}(u), \theta_0,u) \phi_{\theta_0}(u)
\;.
\eeaa
The sums $n^{-1/2} \sum_{i=1}^n U_{1i}(\theta_0)$
and $n^{-1/2} \sum_{i=1}^n U_{2i}(\theta_0)$ are uncorrelated and converge
weakly    to               independent mean zero normal vectors with covariances
$\Sigma_1(\theta_0)$ and $\Sigma_2(\theta_0)$.   
Moreover,
$$
\sqrt n
[\Gamma_{n\hat \theta} - \Gamma_{\theta_0} - [\hat \theta - \theta_0] 
\dot \Gamma_{\theta_0}](t) = {1 \over \sqrt n} \sum_{i=1}^n W_{0i}(t) + o_P(1)
$$
uniformly in $t \in [0, \tau]$.

\item[{(ii )}] We have
$\Sigma_{qn}(\hat \theta) \to_p \Sigma_q(\theta_0)$ for $q = 1,2$,
\beaa
& & \Vert  \Gamma_{n \hat \theta} -   \Gamma_{n \theta_0} \Vert \to_P 0,
\quad \quad  
 \Vert \dot \Gamma_{n \hat \theta} -  \dot \Gamma_{n \theta_0} \Vert \to_P 0
 \;,\\
& & \Vert 
  \int_0^{\cdot} \hat \rho_{\hat \phi}(u, \hat \theta) N_.(du)  
- \int_0^{\cdot} \rho_{\phi_{\theta_0}}(u, \theta_0) EN(du)  \Vert \to_P 0
 \;,\\
& & \Vert 
\int_0^{\cdot} 
{\dot S \over S^2}(\Gamma_{n \hat \theta}(u-), \hat \theta,u) N_.(du) 
- \int_0^{\cdot} {\dot s \over s^2}(\Gamma_{\theta_0}(u-), \theta_0,u) EN(du) 
\Vert \to_P 0 \;,\\
& & \Vert 
\int_0^{\cdot} 
{S' \over S^2}(\Gamma_{n \hat \theta}(u-), \hat \theta,u) N_.(du) 
- \int_0^{\cdot} {s' \over s^2}(\Gamma_{\theta_0}(u-), \theta_0,u) EN(du) 
\Vert \to_P 0 \;,\\
& & \limsup_n \exp  \int_0^{\tau} {|S'| \over S^2}
(\Gamma_{n\hat \theta}(u-), \hat \theta,u) N_.(du) = O_P(1) 
\eeaa
and
$\hat \P_{\hat \theta} (u,t)  \to_P \P_{\theta_0}(u,t)
$ 
uniformly
in $0 < u < t \leq \tau$.

\item[{(iii)}] Let
\beaa
\psi_1(t,D) & = & \pi(D)^{-1}
E 1(Z_i \in D) f(\Gamma_0(t), \theta_0|Z_i) \;,\\ 
\psi_2(t,D) & = & \psi_1(t,X) \dot \Gamma_0(t) + \pi(D)^{-1}
E 1(Z_i \in D) \dot F(\Gamma_0(t), \theta_0|Z_i)  
\eeaa
and let $\hat \psi_p, p = 1,2$ be the estimate of this function
obtained by replacing the pair $(\theta_0, \Gamma_0)$ and the function
$\pi(D)$ by 
$(\hat \theta, \Gamma_{n\hat \theta})$ and $\hat \pi(D)$. Then
$\Vert \hat \psi_q - \psi_q \Vert \to_P 0, q = 1,2$.
\item[{(iv)}] Part (ii) and (iii) remains to hold if the estimates
$(\hat \theta, \Gamma_{n\hat \theta})$ are replaced by
 $(\theta^*, \Gamma^*_n)$ such that $\theta^* \to_P \theta_0$ and
$\Vert \Gamma^*_n - \Gamma_{\theta_0} \Vert \to_P 0$.
\eitem

We omit the proof of this lemma. Part (i)-(ii) and (iv) can be found
in Dabrowska (2005), while part (iii) 
is a straightforward consequence of part (i)-(ii). 

\it Proof of Proposition 1. \rm
We have
$$
\hat W(t,D) = {\pi(D) \over \hat \pi(D)} \sum_{j=1}^4 \hat W_j(t,D) 
\;,
$$
where
\beaa
\hat W_1(t,D) & = & {1 \over \sqrt n} {1 \over \pi(D)}
\sum_{i=1}^n 1(Z_i \in D) [F(\Gamma_0(t),
\theta_0| Z_i) - F_D(t)] \;,\\
\hat W_2(t,D) & = & {1 \over \sqrt n} \sum_{i=1}^n W_{0i}(t)
 \psi_1(t,D) 
   -  \psi_2(t,D)^T \Sigma^{-1}(\theta_0) {1 \over \sqrt n} \sum_{i=1}^n 
U_{2i}(\theta_0) \;,\\
\hat W_3(t,D) & = &  \psi_2(t,D)^T \Sigma^{-1}(\theta_0) {1 \over \sqrt n}
 \sum_{i=1}^n U_{1i}(\theta_0) \\ 
\hat W_4(t,D) & = & {1 \over \sqrt n \pi(D)} \sum_{i=1}^n 1(Z_i \in D)  
[F(\hat \Gamma_{\hat \theta}
(t, \hat \theta|Z_i) - F(\Gamma_{0}
(t), \theta_0|Z_i)] \\
& - & \hat W_2(t,D)  - \hat W_3(t,D) \;.
\eeaa

Here $\hat W_j(t,D), j = 1,2,3$ represent uncorrelated  
sums of mean zero iid processes with finite variance and covariance
\bea
\cov (\hat W_1(t_1,D_1), \hat W_1(t_2,D_2)) & = & {\pi(D_1) \pi(D_2)}^{-1}
E1(Z_i \in D_1 \cap D_2)
F(t_1|Z) F(t_2|Z) \;,  \nonumber\\
& - & F_{D_1}(t_1)F_{D_2}(t_2) \nonumber \\
\cov (\hat W_2(t_1,D_1), \hat W_2(t_2,D_2)) & = &  
\cov (W_0(t), W_0(t')) \psi_1(t_1,D_1)  \psi_1(t_2,D_2)  \nonumber \\
& +  & 
\psi_2(t_1, D_1)^T \cov(W_0(t_1),T) \psi_1(t_2,D_2) \\
& + & 
[\psi_2(t_1, D_1)^T \cov(W_0(t_1),T) \psi_1(t_2,D_2)]^T \nonumber \\
& + & \psi_2(t_1,D_1)^T Var T \psi_2(t_2,D_2) \nonumber \\
& - & 
\cov (\hat W_3(t_1,D_1), \hat W_3(t_2,D_2)) \;, \nonumber \\
\cov (\hat W_3(t_1,D_1), \hat W_3(t_2,D_2)) & = &
\psi_2(t_1,D_1)^T \Sigma^{-1}(\theta_0) \Sigma_1(\theta_0) \Sigma^{-1}
(\theta_0) \psi_2(t_2,D_2) \nonumber 
\eea
and, from section 2, 
\beaa
& & \cov T   =  \Sigma^{-1}(\theta_0),  \quad 
\cov (T, W_0(t))  = 
- \Sigma^{-1}(\theta_0) [\phi_{\theta_0}+ \dot \Gamma_{\theta_0}](t) 
\;,\\ 
& & \cov(W_0(t), W_0(t'))   =   K_{\theta_0}(t,t') \;.
\eeaa
We also have 
$$
[\phi_{\theta_0}+ \dot \Gamma_{\theta_0}](t) =
\int_0^{\tau} K_{\theta_0}(t,u) \rho_{\phi}(u, \theta_0) EN.(du)
\;.
$$
By central limit theorem, finite dimensional distributions of the processes
$\{\hat W_j, j = 1,2,3\}$ converge weakly to 
a multivariate vector with covariance matrix given by 
(14).

For each $j = 1,2,3$, the process $\{\hat W_j(t,D): (t,D) \in \T\}$
 can be represented
as  \newline
 $n^{-1/2} \sum_{i=1}^n h^{(j)}_{t,D}(X_i, \delta_i, Z_i)$, 
with $h^{(j)}$ varying over a  
Euclidean class of functions
$\H_j = \{h_{t,D}^{(j)}: (t, D)  \in \T\}$ for a square integrable envelope
(Nolan and Pollard, 1987).
This can be verified, by noting that $\D$ is a finite collection
of sets, and for each $D \in \D$, the relevant functions 
$h_{t,D}^{(j)} \in \H_j$ can be represented
as finite  linear 
combination of  functions of bounded variation with respect to t.
We also have $Eh_{t,D}^{(j)} (X_i, \delta_i, Z_i) = 0$ for each 
$h_{t,D}^{(j)} \in \H_j$.
Hence the process $\hat W_j = G_{n,j} = \{ \sqrt n [P_n - P](h_{t,D}^{(j)}):  
h_{t,D}^{(j)} \in \H_j\}$ 
is  equicontinuous and $\H_j$ is totally bounded 
with respect to the variance semi-metric $\rho_j$. Set $\rho = \max {\rho_j,
j =1,2,3}$. Then $\T$ is totally bounded with respect to $\rho$, and
$\{ \hat W_j: j = 1,2,3\}$ is asymptotically tight in 
$\l^{\infty}(\T)$ and
 converges weakly to a
Gaussian process $\{W_j: j =1 ,2, 3\}$.
Its components are independent, and $W_j$ have covariance function
given by the right--hand side of (14).

Using Taylor expansion, we also have
$\hat W_4(t,D) = \hat W_{41}(t,D) + \hat W_{42}(t,D)$, where
\beaa
\hat W_{41}(t,D) & = & \sqrt n(\Gamma_{n\theta}(t) - \Gamma_{\theta_0}(t) - 
(\hat \theta - \theta_0)^T \dot \Gamma_0(t))  \hat \psi^*_1(t, D)  \\
& - & {1 \over \sqrt n} \sum_{i=1}^n W_{0i}(t) 
\psi_1(t,D) \;,\\
\hat W_{42}(t,D) & = & \psi_2^*(t,D)^T \sqrt n(\hat \theta - \theta_0) 
- \psi_2(t,D)^T \Sigma^{-1}(\theta_0) \sqrt n U_n(\theta_0) 
\eeaa
and
\beaa
\psi_1^*(t,D) & = & {1 \over n \pi(D)} 
\sum_{i=1}^n 1(Z_i \in D) f(\Gamma^*(t), 
\theta^*|Z_i) \;,\\
\psi_2^*(t,D) & = &  \psi_1^*(t,D) \dot \Gamma_0(t) + 
{1 \over n \pi(D)} \sum_{i=1}^n 1(Z_i \in D) \dot F(\Gamma^*(t), 
\theta^*|Z_i) \;.
\eeaa
Here $\theta^*$ is on a line segment between $\theta_0$ and
$\hat \theta$, and $\Vert \Gamma^* - \Gamma_{\theta_0} \Vert \to_P 0$.  
By  Lemma 1, 
$$
\sup \{|\hat W_{4p}(t,D)|: (t,D) \in \T\} \to_P 0
$$
for $p = 1,2$. To complete the proof of part (i) of the Proposition 3, 
we note that 
$\hat \pi(D) \to \pi(D)$ a.s. for $D \in \D$ so that 
$\hat W = \{(\pi(D) /\hat \pi(D)) \sum_{j=1}^4 \hat W_j(t,D): (t,D) \in \T\}$
converges weakly in $\l^{\infty}(\T)$ to  $W = \{
W(t,D) =  \sum_{j=1}^3 W_j(t,D): (t,D) \in \T\}$. 
Its variance function
is given by 
$v_D(t)  =  \sum_{j=1}^3 \var W_j(t,D)$. For any $D$ this is a 
continuous function with respect to $t$ and positive on any interval
$[\tau_1 - \epsilon, \tau_2 + \epsilon]$ on which 
$\Gamma_{\theta_0}
$ forms
a continuous strictly increasing function.

To show part (ii), 
first recall that $V_i = (V_{1i}, 
V_{2i}), i = 1, \ldots,n, \ldots$ and $V_{3} = (V_{31}, \ldots, V_{3d})$,  
are mutually independent 
$\N(0,1)$ variables, independent of 
$R_i = (X_i, \delta_i,Z_i),
i = 1, \ldots, n$. 
We let variables $R_i =  i = 1,2, \ldots$ be defined
as coordinate projections on the ``first'' $\infty$ coordinates in the 
product probability space $(\Omega^{\infty} \times \V
\times \V', \F^{\infty} \times
\B \times B', P^{\infty} \times Q \times Q')$ and 
let $V_i, i = 1, \ldots,..$ and $V_3$ be defined
on the ``last'' two coordinates.

Set
\beaa
\tilde  W_1(t,D) & = & {1 \over \sqrt n} {1 \over \pi(D)}\sum_{i=1}^n V_{1i} 1(Z_i \in D)
[F(\Gamma_{\theta_0}(t),  \theta_0|Z_i) - F_D(t)] \;,\\
\tilde W_2(t,D) & = & 
\tilde W_0(t)  \psi_1(t,D) + \int_0^{\tau}
\tilde  W_0(s)  \rho_{\phi_{\theta_0}}(s,  \theta_0) EN_.(ds) 
\Sigma^{-1}(\theta_0)  \psi_2(t,D) \;,\\
\tilde  W_3(t,D) & = &  V_3  \Sigma_{1}^{1/2} (\theta_0) 
\Sigma^{-1}(\theta_0) \psi_1(t,D) \;,
\eeaa
where 
$$
\tilde W_0(t) = {1 \over \sqrt n} \sum_{i=1}^n V_{2i} {1[X_i \leq t, 
\delta_i = 1] \over s(\Gamma_{\theta_0}(X_i-),  \theta_0, X_i)} 
\P_{\theta_0}(X_i, t) \;.
$$

For $j,k  = 1,2,3, j \not=k$, we have
\beaa
 & & \cov (\tilde W_j(t,D), \tilde W_j(t',D')) =
 \cov (\hat W_j(t,D),  \hat W_j(t',D')) \;,\\
 & & \cov (\hat W_k(t,D), \tilde W_j(t',D')) = 
 \cov (\tilde  W_k(t,D), \tilde W_j(t',D')) = 
 \cov (\hat W_k(t,D), \hat W_j(t',D')) = 0 \;.   
\eeaa
Also $\tilde W_3$ does not involve $n$, the $R_i, i = 1,2, \ldots$ 
or the  $V_{ji}, j = 1,2, i = 1,2 \ldots$ sequences, and is 
independent of the processes $\tilde W_j, j =1,2$ and $\hat W_j, j = 1,2,3$.

Similarly to part (i),
the processes $\{\tilde W_j(t,D): (t,D) \in \T,j = 1,2\}$ 
are of the form
$\tilde W_j(t,D) = {1 \over \sqrt n} \sum_{i=1}^n V_{ji} 
g^{(j)}(X_i, \delta_i,Z_i)$,
where $g^{(j)}$ varies over  
$\G_j = \{g_{t,D}^{(j)}(x,d,z): (t,D) \in \T\}$, a Euclidean class of functions
for a square integrable envelope and  is totally bounded with respect 
to the semi-metric $\rho$. The class of products 
$\{vg_{t,D}^{(j)}(x,\delta,z): (t,D) \in \T\}$ 
is also Euclidean. Therefore, unconditionally 
$[\tilde W_j: j = 1,2]$ is asymptotically tight and 
converges to a Gaussian process $[W_j^{\#}: j = 1,2]$, whose
components are independent and independent of $\tilde W_3$
and $[W_1, W_2,W_3]$.

Alternatively,
for $j=1$, we have $g^{(1)}_{t,D} = 
h^{(1)}_{t,D}$ with $Ph^{(1)}_{t,D} = 0$ and
$$
\tilde W_1(t,D) = 
{1 \over \sqrt n}\sum_{i=1}^n V_{1i}(\delta_{R_i} - P)[g_{t,D}] 
= {1 \over \sqrt n}\sum_{i=1}^n V_{1i}\delta_{R_i}[g_{t,D}] \;.
$$
For $j = 2$
\beaa
\tilde W_2(t,D) & = & {1 \over \sqrt n}\sum_{i=1}^n V_{2i}(\delta_{R_i} - P)
[g^{(2)}_{t,D}] + {1 \over \sqrt n} \sum_{i=1}^n V_{2i} P 
[g^{(2)}_{t,D}] \\
& = & \tilde W_{21}(t,D) + \tilde W_{22}(t,D) 
\eeaa
and the two components on the right-hand side are uncorrelated.
Application of the unconditional
multiplier central limit theorem in van der Vaart and Wellner (1996,
Corollary 2.9.4, p.180) implies that the processes 
$[\tilde W_1, \tilde W_{21}, \tilde W_{22}, \tilde W_3]$ and
$[\hat W_1, \hat W_{2}, \hat W_3]$ converge jointly 
 in $[\l^{\infty}(\T)]^4 \times [\l^{\infty}(\T)]^3$ 
to independent Gaussian processes, $[W_1^{\#}, W_{21}^{\#},
W_{22}^{\#},  W_3^{\#} = \tilde W_3]$ and $[W_1, W_2,W_3]$.
By continuous mapping theorem, we also have unconditional
weak convergence of $[\hat W = \sum_{j=1}^3 \hat W_j, 
\tilde W = \sum_{j=1}^3 \tilde W_j]$ 
in $\l^{\infty}(\T)  
\times \l^{\infty}(\T)$ to a vector of independent Gaussian
processes $[W, W^{\#}]$, with the same covariance function.

Conditionally on $R_1, R_2, \ldots,...$ the processes $\tilde W_1, \tilde 
W_{21}$ and $\tilde W_{22}$ have mean zero,  
\beaa
\cov_V [\tilde W_1(t_1,D_1), \tilde W_1(t_2,D_2)] & = & {1 \over n} \sum_{i=1}^n
g_{t_1, D_1}^{(1)}(R_i)g^{(1)}_{t_2, D_2}(R_i)^T \to 
P(g_{t_1, D_1}^{(1)}[g^{(1)}_{t_2, D_2}]^T) 
\;, \\ 
\cov_V [\tilde W_{21}(t_1,D_1), \tilde W_{21}(t_2,D_2)] & = & 
{1 \over n} \sum_{i=1}^n
g_{t_1, D_1}^{(2)}(R_i)g_{t_2, D_2}^{(2)}(R_i)^T - Pg_{t_1, D_1}^{(2)} 
[Pg_{t_2, D_2}^{(2)}]^T \\ 
&\to & \cov(
g_{t_1, D_1}^{(2)}(R_1),g_{t_2, D_2}^{(2)}(R_1)) \;,\\
\cov_V [\tilde W_{22}(t_1,D_1), \tilde W_{22}(t_2,D_2)] & = & 
 Pg_{t_1, D_1}^{(2)} 
[Pg_{t_2, D_2}^{(2)}]^T \;, \\ 
\cov_V [\tilde W_{21}(t_1,D_1), \tilde W_{22}(t_2,D_2)] & = &  
{1 \over n} \sum_{i=1}^n[
g_{t_1, D_1}^{(2)}(R_i) - Pg_{t_1, D_1}^{(2)}] 
[Pg_{t_2, D_2}^{(2)}]^T \to 0 \;,\\ 
\cov_V [\tilde W_{1}(t_1,D_1), \tilde W_{2j}(t_2,D_2)] & = & 0, 
\quad j=1,2 \;, \\
\cov_V [\tilde W_{3}(t_1,D_1), \tilde W_{2j}(t_2,D_2)] & = & 0, 
\quad j=1,2 \;,\\
\cov_V [\tilde W_{3}(t_1,D_1), \tilde W_{1}(t_2,D_2)] & = & 0, 
\eeaa
for almost all $R_1, R_2, \ldots $. 
(Actually, conditionally on $R_1, R_2, \ldots $, $\tilde W_j$ processes
are independent).
By conditional multiplier CLT,
we have that 
conditionally on $R_1, R_2, \ldots,$ the finite dimensional
distributions of $\tilde W_1$ and  $\tilde W_{2}$
are asymptotically multivariate normal and independent, for almost all 
$R_1, R_2 \ldots$. The covariance function is the same as of finite
dimensional distributions of $W_1$ and $W_2$. By continuous mapping theorem,
we also have that conditionally on $R_1, R_2, \ldots$, the finite dimensional
distributions of 
$\tilde W$ converge weakly to a multivariate normal
distribution for almost all $R_1, R_2, \ldots$. The covariance
of the multivariate normal distributions is the same as the 
covariance of the corresponding finite dimensional distributions
of $W$.

Let $BL_1$ be the collection of functions $f$ from $\l^{\infty}(\T)$
into $[0,1]$ that are Lipschitz continuous with Lipschitz continuity
constant equal to 1. For fixed $\delta$ 
and $x \in \T$, let $\Pi_{\delta}(x)$ be the 
closest point to $x$  in $\T$ in a partition of the set $\T$ with mesh-width
$\delta$ (with respect to the semi-metric $\rho$).  
By triangular inequality
\beaa
& & \sup_{f \in BL_1} |E_V f(\tilde W) - E f(W)| \leq 
\sup_{f \in BL_1} |E f(W \circ \Pi_{\delta} ) - E f(W)| + \\
& & \sup_{f \in BL_1} |E f(W \circ \Pi_{\delta}) - 
E f_V(\tilde W \circ \Pi_{\delta})| + 
\sup_{f \in BL_1} |E_V f(\tilde W \circ \Pi_{\delta}) - 
E_Vf (\tilde W)| \\
& & = I_1 + I_2 + I_3  \;.
\eeaa
As in van der Vaart and Wellner (1996, p. 182),
the term $I_1$ converges to 0, because the process $W$ has continuous paths
with respect $\rho$ and 
$W \circ \Pi_{\delta} \to W$ in almost surely as $\delta \downarrow 0$.
For fixed $\delta > 0$, $I_2$ converges to 0 for almost all
$R_1, R_2, \ldots$. This follows because conditionally on
$R_1, R_2, \ldots$, the finite dimensional
distributions of $\tilde W$ converge in distribution to 
a multivariate normal vector, for almost all $R_1, R_2, \ldots$.
Finally, 
\beaa
I_3 & \leq & 
\sup_{f \in BL_1} E_V |f(\tilde W \circ \Pi_{\delta}) - 
f (\tilde W)| \leq E_{V_1} 
\Vert \tilde W_1 \circ \Pi_{\delta} - \tilde W_1 \Vert_{\G_{1\delta}} + \\
& + & E_{V_2}\Vert \tilde W_{2} \circ \Pi_{\delta} - \tilde W_{2} \Vert_{\G_{2\delta}}
+ E_{V_3} \Vert \tilde W_{3} \circ \Pi_{\delta} - \tilde W_{3} \Vert_{\G_{3\delta}} \\
& \leq &
E_{V_1} \Vert \tilde W_1 \Vert_{\G_{1\delta}}
+ E_{V_2} \Vert \tilde W_{2} \Vert_{\G_{2\delta}} +
E_{V_3} \Vert \tilde W_{3} \Vert_{\G_{3\delta}} \;,
\eeaa
where  $\G_{j\delta} = \{g - g': g,g' \in \G_j: \rho(g - g') < \delta\}$, 
for $j = 1,2,3$.
The first two expectation converge to 0 as $n \to \infty$ and 
$\delta \downarrow 0$, 
by Lemma 2.9.1
in van der Vaart and Wellner (1996, p 177). 
The last expected does not depend on n, and converges to 0 as 
$\delta \downarrow 0$.

It remains to consider the process 
$\hat W^{\#}$ defined in Section 2.
We  show that unconditionally 
$\Vert \hat W^{\#}_j - \tilde W_j \Vert \to 0$ in probability.
If this is the case, then for $\epsilon > 0$, we have
\beaa
\sup_{f \in BL_1} |E_V^*f(\hat W^{\#}) - Ef(W) | & \leq &
 \sup_{f \in BL_1} |E_V f(\tilde  W) - Ef(W) | 
 +  \sup_{f \in BL_1} |E_V^* f(\hat W^{\#}) - E_V f(\tilde W) | \\
&\leq &   \sup_{f \in BL_1} |E_V f(\tilde  W) - Ef(W) | 
+ \epsilon + 2P_V^*( \Vert \hat W^{\#} - \tilde W \Vert > \epsilon) \;.
\eeaa
The first term converges to 0 in  probability. 
The last term converges to 0 in (outer)
mean.

Clearly, for $j = 3$, we have
 $\hat \Sigma_n(\hat \theta) \to
\Sigma(\theta_0)$, $\hat \Sigma_{2n}(\hat \theta) \to 
\Sigma_2(\theta_0)$ and $\Vert \hat \psi_1 -  \psi_1 \Verti \to 0$
in probability so that $\Vert \hat W^{\#}_3 - \tilde W_3 \Vert \to_P 0$.

Next, for $j = 1,2,3$, define 
$$
\tilde H_j(t,D) = {1 \over n} \sum_{i=1}^n V_{1i} 1(Z_i \in D) h_{jt}(Z_i)
\;,
$$
where 
\beaa
h_{jt}(Z) & = & \pi(D)^{-1} \quad \quad \quad \quad \quad \quad \quad j=1 \;, \\
& = & \pi(D)^{-1} f(\Gamma_{\theta_0}(t), \theta_0|Z) \; \quad j=2 \;,\\
& = & \pi(D)^{-1} \dot F(\Gamma_{\theta_0}(t), \theta_0 | Z) \quad j=3 \;. 
\eeaa
We have $E\tilde H_j(t,D) = 0$ for $(t,D) \in \T$. Unconditionally, 
the strong law of large numbers, yields $\tilde H_j(t,D) \to 0$ a.s.
pointwise in $(t,D) \in \T$. The convergence is also uniform 
since for each $D$, the process $H_j(t,D)$ has paths of bounded variation.
We also have
  $\tilde W_1 - \hat W^{\#}_1 = \sum_{j=1}^4 \tilde W_{1j}$, where
\beaa
\tilde W_{11}(t,D) & = & - \sqrt n[\hat F_D - F_D](t) \tilde H_1(t,D) \;,\\
\tilde W_{12}(t,D) & = & \sqrt n[\Gamma_{n\hat \theta} - \Gamma_0 
- (\hat \theta - \theta_0)^T \dot \Gamma_{\theta_0}](t)\tilde H_2(t,D) \;,\\
\tilde W_{13}(t,D) & = & \sqrt n[\hat \theta - \theta_0]^T  
[\dot \Gamma_{\theta_0}(t) \tilde H_2(t,D) +
 \tilde H_3(t,D)] \;,\\ 
\tilde  W_{14}(t,D) & = & O_P(1) {1 \over n} \sum_{i=1}^n |V_{1i}|
O(\sqrt n \Vert \Gamma_{n\hat \theta}-\Gamma_{\theta_0} \Vert^2 + 
\sqrt n (\hat \theta - \theta_0)^2) \;.
\eeaa
These four terms satisfy $\Vert \tilde W_{1j} \Vert \to 0$ in probability 
(unconditionally)
 and the same holds for the process 
$\tilde W_1 - \hat W_1^{\#}$.

Finally, define 
\beaa
\tilde M(t) & = & {1 \over \sqrt n} 
\sum_{i=1}^n V_{2i} 1(X_i \leq t, \delta_i =1) \;,\\
\tilde W_4(t) & = & {1 \over \sqrt n} \sum_{i=1}^n V_{2i} {
1(X_i \leq t, \delta_i=1) \over s(\Gamma_0(X_i-), \theta_0, X_i)} = \int_0^t
{\tilde M(du) \over s(\Gamma_{\theta_0}(u-), \theta_0,u)} \;,\\
\hat W_4^{\#}(t) & = & {1 \over \sqrt n} \sum_{i=1}^n V_{2i} {
1(X_i \leq t, \delta_i=1) \over S(\Gamma_{n\hat \theta}(X_i-), \hat \theta, X_i)} = \int_0^t
{\tilde M(du) \over S(\Gamma_{n\hat \theta}(u-), \hat \theta,u)} \;.
\eeaa
A similar argument as in analysis of the term $\tilde W_2$ shows that
$\tilde W_4$ converges weakly (unconditionally) to a mean zero time transformed
Brownian motion with variance function $C_{\theta_0}(t)$. Since $EN$ is
a continuous function,   so is $C_{\theta_0}$. We have
$$
\hat W_4^{\#}(t) - \tilde W_4(t) = \int_0^t \left 
[{s(\Gamma_{\theta_0}(u-), \theta_0,u)
\over S(\Gamma_{n\hat \theta}(u-), \hat \theta,u)} - 1 \right] 
\tilde W_4(du) \;.
$$
Denote the term in the bracket by $\hat a_n(u-)$. Then $\hat a_n$ is a process
with left continuous and right-hand limits, 
$\Vert a_n \Vert \to_P 0$ and
$$
\limsup_n \Vert a_n \Vert_v  = O_P(1) \;,
$$
where $\Vert \cdot \Vert_v$ is 
the variation norm. For given $\delta   > 0$, let
$t_1 < t_2 < \ldots t_k$ be a partition of $[0, \tau]$, such 
that $C_{\theta_0}(t_i)
- C_{\theta_0}(t_{i-1}) < \delta$. Define $\Pi_{\delta}(t) = t_{i-1}$
if $t \in [t_{i-1}, t_i)$.
Then integration by parts, yields
\beaa
\hat W_4^{\#}(t) &= &
\int_0^t a_n(u-)[
\tilde W_4 - \tilde W_4 \circ \Pi_{\delta}](du) +
\int_0^t a_n(u-) [\tilde W_4 \circ \Pi_{\delta}](du) \\
& = & [\tilde W_4 - \tilde W_4 \circ \Pi_{\delta}](t) a_n(t) + \int_0^t
[\tilde W_4 - \tilde W_4 \circ \Pi_{\delta}](u) a_n(du)+
\int_0^t a_n(u-)  [\tilde W_4 \circ \Pi_{\delta}](du) \;.
\eeaa
The right-hand side converges then to 0 in probability
uniformly in $t$,
as $n \to \infty$, followed by $\delta \to 0$. 
We also have
\beaa
\hat W^{\#}_0(t) & = & \int_0^t \hat W^{\#}_4(du) \P_{\hat \theta}(u,t) 
\;,\\
\tilde W_0(t) & = & \int_0^t \tilde W_4(du) \P_{\theta_0}(u,t) \;.
\eeaa
Then
\beaa
\hat W^{\#}_0(t) & = & \hat W_4^{\#}(t) - \int_0^t \hat W^{\#}_0(u-) 
{S' \over S^2}(\Gamma_{n\hat \theta}(u-), \hat \theta,u) N_.(du) \;,\\
\tilde W_0(t) & = & \tilde W_4(t) - \int_0^t \tilde W_0(u-) 
{s' \over s^2}(\Gamma_{\theta_0}(u-),  \theta_0,u) EN(du) \;.
\eeaa
We have
\beaa
[\hat W^{\#}_0(t) - \tilde W_0(t) ]  & = & \mbox{Rem}(t) -
 \int_0^t [\hat W^{\#}_0 - \tilde W_0](u-) 
{S' \over S^2}(\Gamma_{n \hat \theta}(u-), \hat \theta,u) N_.(du) \;, \\
\mbox{Rem}(t)  & = &  
[\hat W^{\#}_4 - \tilde W_4](t) \\
& - & \int_0^t \tilde W_0(u-) \left (
{S' \over S^2}(\Gamma_{n \hat \theta}(u-), \hat \theta,u) N_.(du) 
- {s' \over s^2}(\Gamma_0(u-), \theta_0,u) EN(du) \right ) \;.
\eeaa

We have
$\Vert \mbox{Rem} 
\Vert \to 0$ and $\Vert \mbox{Rem}^- \Vert \to 0$  in probability.
Hence by Gronwall's inequality (Beesack (1975))
\beaa
& & |\hat W^{\#}_0 - \tilde W_0|(t)   \leq   |\mbox{Rem}(t)| \\
& + &  \int_0^t |\mbox{Rem}(u-)| {|S'| \over S^2}
(\Gamma_{n\hat \theta}(u-), \hat \theta,u) N_.(du) 
\exp \int_u^{\tau} {|S'| \over S^2}
(\Gamma_{n\hat \theta}(u-), \hat \theta,u) N_.(du) \\
& \leq &\max {\sup_{t \leq \tau} |\mbox{Rem}(t)|, |
\mbox{Rem}(t-)|}
\limsup_n  \exp \int_0^{\tau} {|S'| \over S^2}
(\Gamma_{n\hat \theta}(u-), \hat \theta,u) N_.(du)      
\;.
\eeaa
Application of Lemma 1 and integration by parts
implies that this term converges
to 0 in probability, and 
 $\Vert \hat W_0^{\#} - \tilde W_0 \Vert \to 0$ in probability.
 Similarly, we have  
$\Vert \hat W^{\#}_2 - \tilde W_2 \Vert \to 0$ in probability. \done

\noindent
\bf Acknowledgement. \rm I thank an anonymous reviewer and 
Roger Koenker for comments. 

\bref

\item[
Andersen, P. K. and Gill, R. D.] (1982). Cox's regression model
for counting processes: a large sample study. \it Ann. Statist.
\bf 10 \rm 1100-1120.

\item[
Bennett. S.] (1983). Analysis of the survival data by the proportional
odds model. \it Statistics in Medicine, \bf 2 \rm 273--277.

\item[
Beesack, P. R.] (1975). \it Gronwall Inequalities. \rm Carlton Math. Lecture
Notes \bf 11, \rm Carlton University, Ottawa.

\item[
Bie, O., Borgan, O. and Liestol, K.] (1987). Confidence intervals and 
confidence bands for the cumulative hazard rate function and their 
small-sample properties. \it Scand. J. Statist. \bf 14 \rm 221--233.

\item[
Burr, D. and Doss, H.] (1993). Confidence bands for the median survival
time as a function of the covariates in the Cox model. \it J. Amer.
Statist. Assoc. \bf 88 \rm 1330--1340.

\item[
Bogdanovicius, V. and Nikulin, M.] (1999). 
Generalized proportional hazards model based on modified partial likelihood.
\it Lifetime Data Analysis \bf 5 \rm 329-350.

\item[
Cheng, S. C., Wei, L. J. and Ying, Z.] (1995) Analysis of transformation
models with censored data. \it Biometrika, \bf 82 \rm 835-845.

\item[Cox, D. R.] (1972). Regression models in life tables. \it J. Roy.
Statist. Soc. Ser. B. \bf 34 \rm 187--202.

\item[
Dabrowska, D. M. and Doksum, K. A.] (1987). Estimates and confidence intervals
for median and mean life in the proportional hazard model. \it
Biometrika, \bf 74, \rm 799-807.

\item[
Dabrowska, D. M.] (2005) Estimation in a class of 
semiparametric transformation models.  \it J. Multivariate. Analysis. \rm
(in revision).

\item[
Hjort, N.] (1985) Bootstrapping Cox's regression model. Technical
Report 241. Stanford University, Dept Statistics.

\item[ 
Kalbfleisch, J. D. and Prentice, R. L.] (1980). \it The Statistical Analysis
of Failure Time Data. \rm New York, Wiley

\item[
Koenker, R. and Geling, O.] (2001). Reappraising medfly
longevity: a quantile regression survival analysis. \it
J. Amer Statist. Assoc. \bf 96 \rm 458-468.

\item[
Lin, D. Y., Fleming, T. R., Wei, L.J.] (1994). Confidence bands for
survival curves under the proportional hazards model. \it Biometrika,
\bf 81 \rm 73-81.

\item[
Murphy, S. A., Rossini, A. J. and  van der Vaart, A. W.] (1997). Maximum
likelihood estimation in the proportional odds model. \it
J. Amer. Statist. Assoc. \bf 92, \rm 968--976.

\item[
Nolan, D. and Pollard, D.] (1987). U-processes: rates of convergence.
\it Ann. Statist. \bf 15 \rm 780-799.

\item[
Pettitt, A.N.] (1984). Proportional odds models for survival
data and estimates using ranks. \it Applied Statistics \bf 33 \rm 169-175.

\item[Portnoy, S.] (2003). Censored regression quantiles. \it J. Amer.
Statist. Assoc. \bf 98 \rm 1001--1013.

\item[
van der Vaart, A.W. and Wellner, J.A.] (1996). \it Weak convergence
and Empirical Processes with Applications to Statistics. \rm Springer Verlag
\eref

\newpage

\begin{center}
Table 1. Regression estimates and standard errors \newline
in the proportional odds ratio model. 
\end{center}

\begin{table}[ht]
\centerline{
\begin{tabular}{ccccccc} 
\hline 
\hline
& & & & \\
covariate & & theta & sd error & p-value \\
\hline
\hline
& & & &    \\
PS        & & -1.049 &     0.045 &           $< 10^{-5}$      \\
 SQUAMOUS & &  -0.246 &    0.428 &         0.71      \\
 SMALL    & &   1.345 &     0.304 &         0.01      \\
 ADENO    & &   1.275  &    0.342  &        0.02      \\
 LARGE    & &     NA    &     NA     &           NA  \\
\hline
\hline 
\end{tabular} }
\end{table}

\newpage

\noindent
Figure captions:

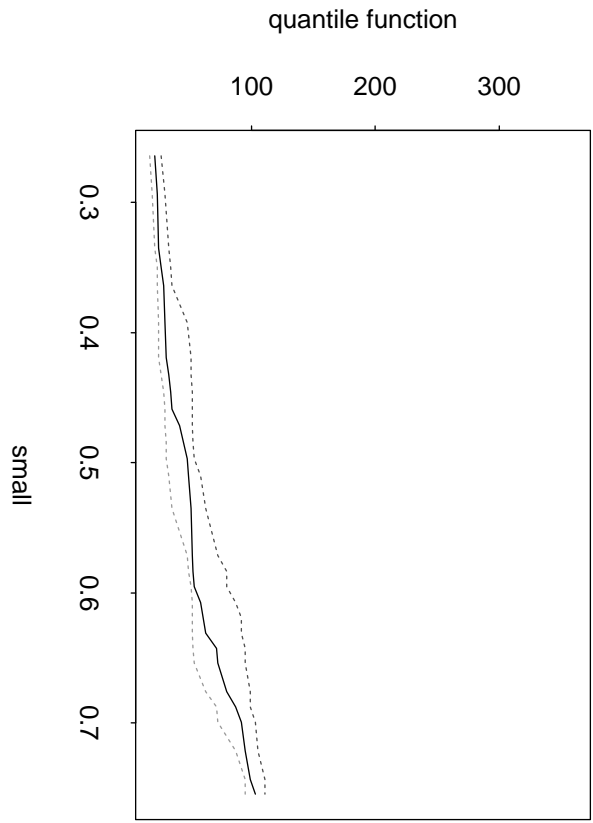
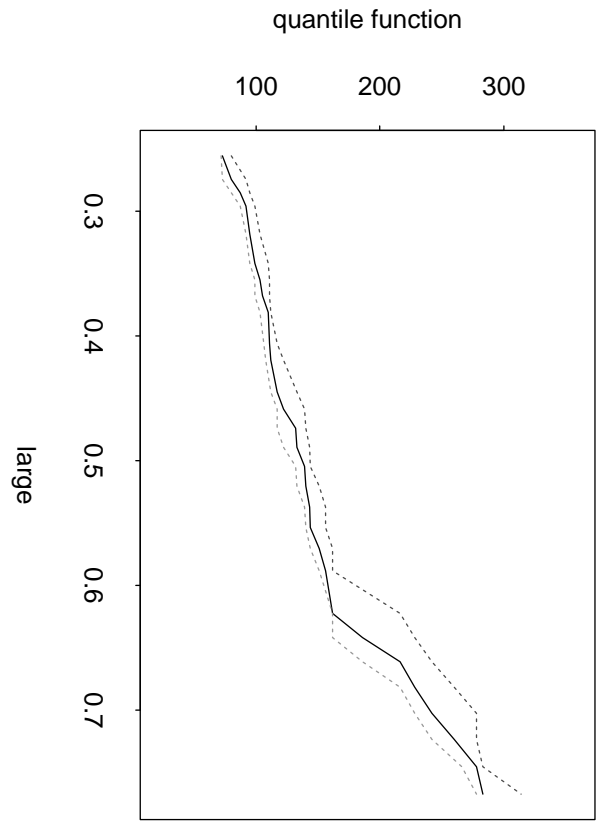
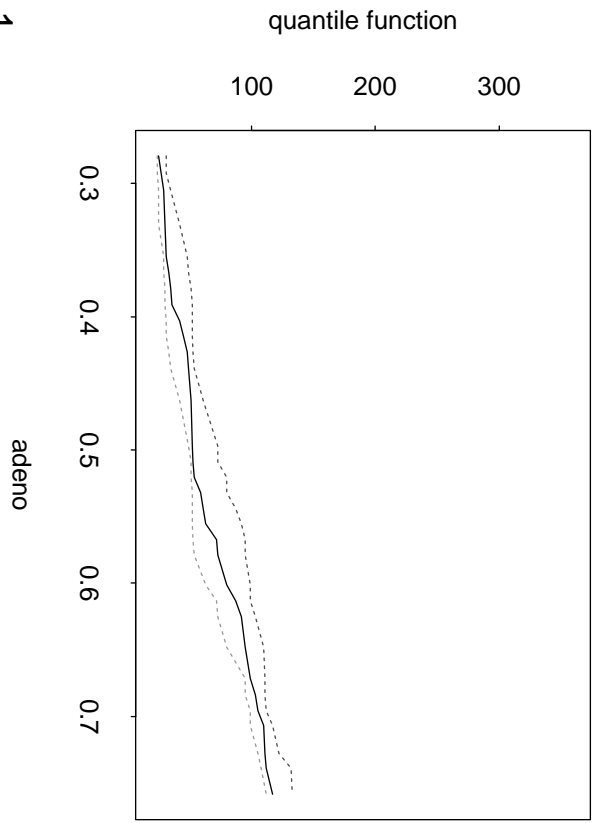
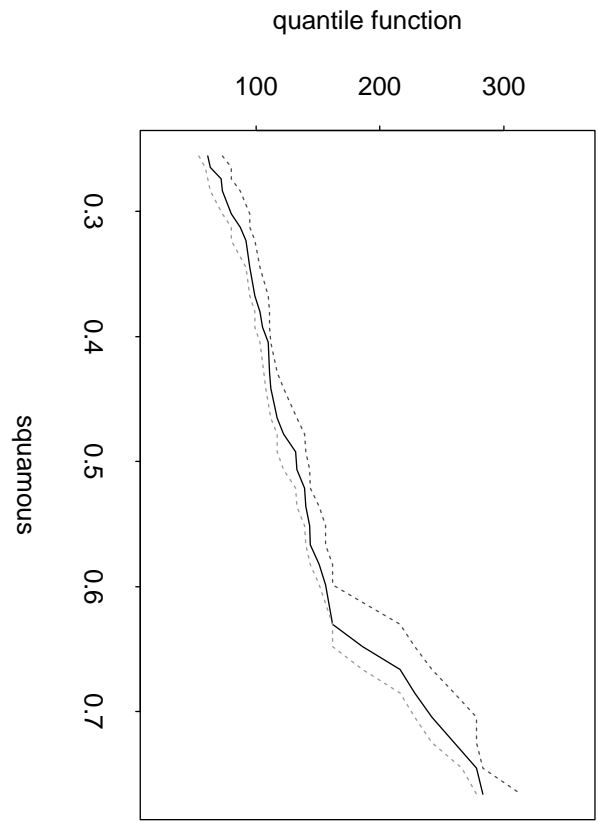
Figure 1. Quantile regression and 95\% simultaneous confidence 
bands. Covariate space partitioned according to four tumor types.

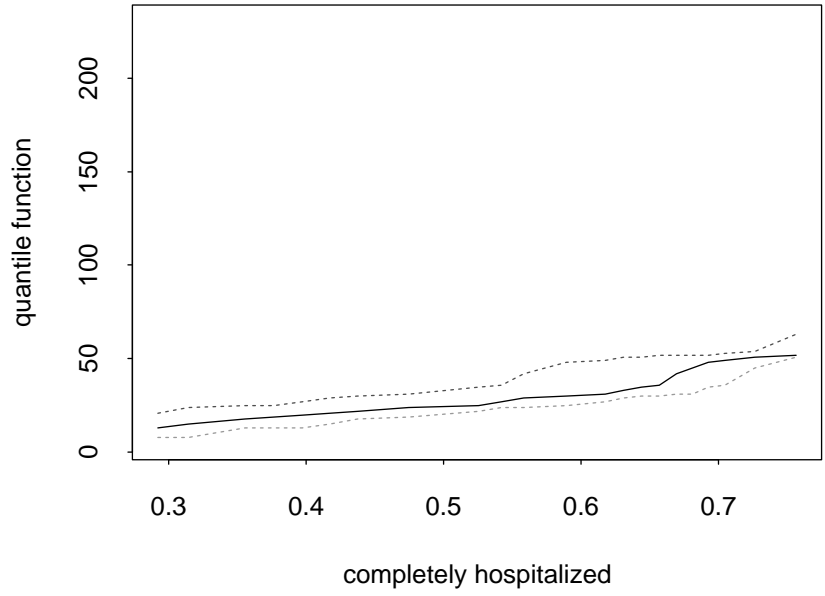
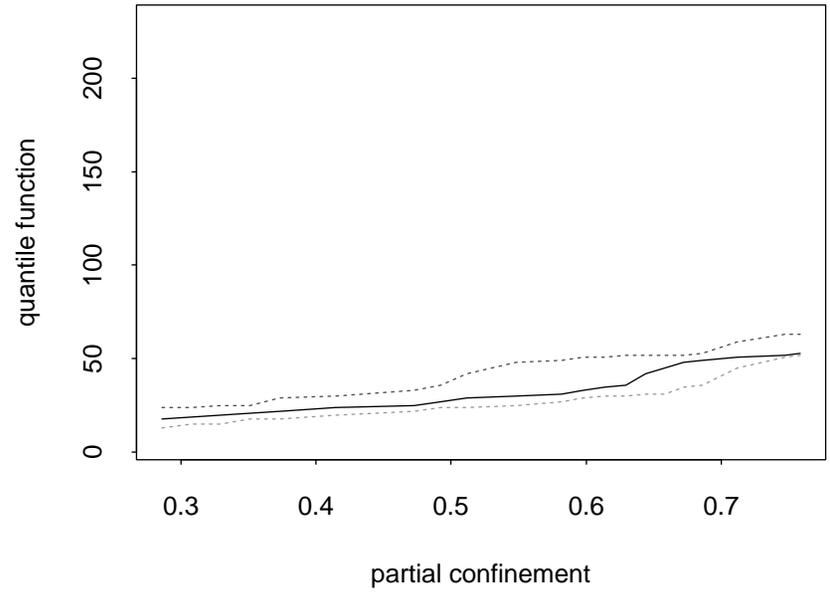
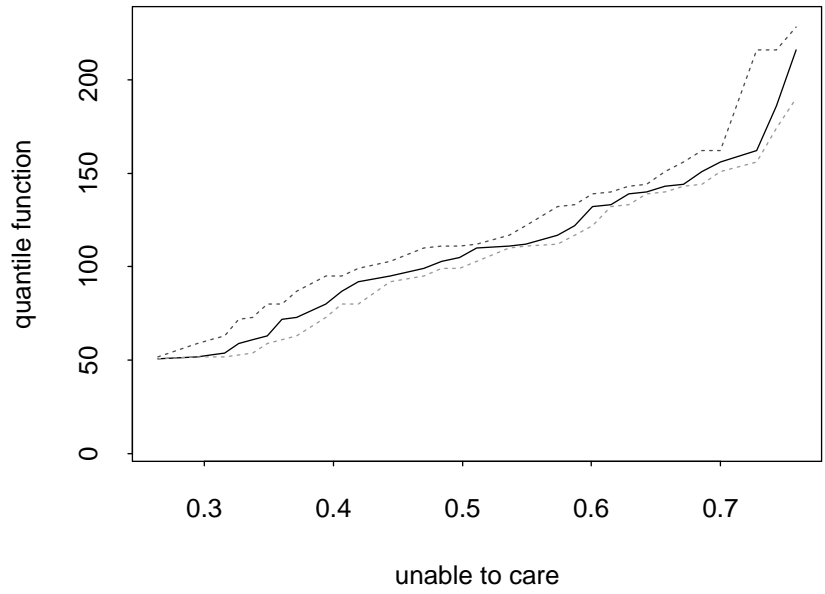
Figure 2. Quantile regression and 95\% simultaneous confidence 
bands. Covariate space partitioned into three groups
according to the 
of performance status (Karnofsky) score. 

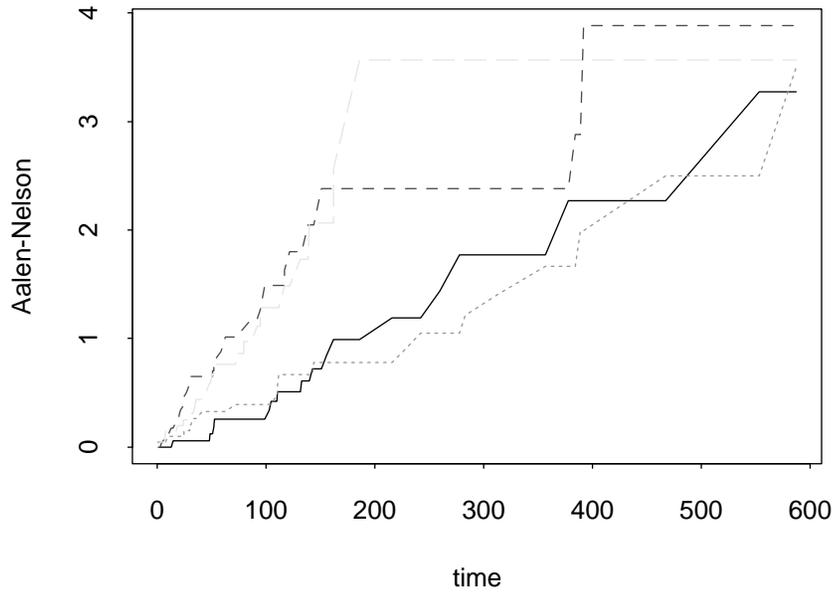
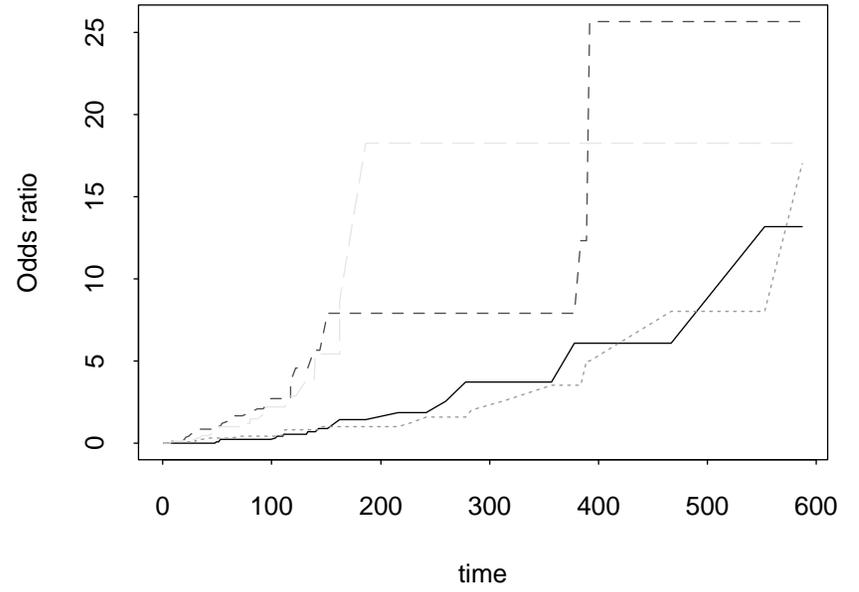
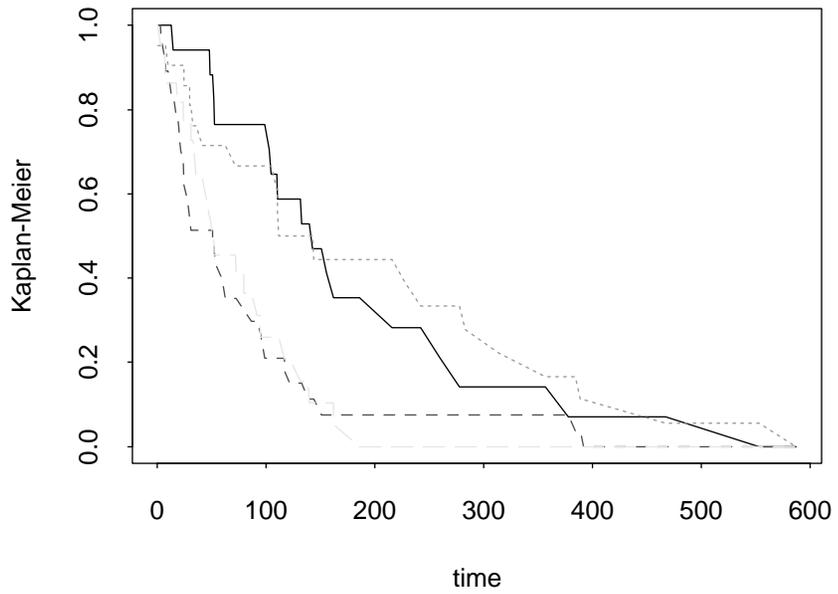
Figure 3. Aalen-Nelson,  odds ratio
function and Kaplan-Meier estimators for the four
tumor cell types: squamous  (solid line), large (dotted line),
small  (long dash) and adeno (short dash).

\end{document}